\documentclass[12pt]{amsart}
\usepackage{amssymb}
\newtheorem{theorem}{Theorem}[section]
\newtheorem{lemma}[theorem]{Lemma}
\newtheorem{proposition}[theorem]{Proposition}
\newtheorem{example}[theorem]{Example}
\newtheorem{definition}[theorem]{Definition}
\newtheorem{remark}[theorem]{Remark}
\newtheorem{corollary}[theorem]{Corollary}

\newcommand{\B}{\mbox{$\mathbb{B}$}}
\newcommand{\C}{\mbox{$\mathbb{C}$}}
\newcommand{\Kp}{\mbox{$\mathbb{K}$}}
\newcommand{\M}{\mbox{$\mathbb{M}$}}
\newcommand{\N}{\mbox{$\mathbb{N}$}}
\newcommand{\Tr}{\mbox{$\mathbb{T}$}}

\newcommand{\A}{\mbox{${\mathcal A}$}}
\newcommand{\Br}{\mbox{${\mathcal B}$}}
\newcommand{\CB}{\mbox{${\mathcal{CB}}$}}
\newcommand{\E}{\mbox{${\mathcal E}$}}
\newcommand{\F}{\mbox{${\mathcal F}$}}
\newcommand{\Hi}{\mbox{${\mathcal H}$}}
\newcommand{\K}{\mbox{${\mathcal K}$}}
\newcommand{\Li}{\mbox{${\mathcal L}$}}
\newcommand{\LM}{\mbox{${\mathcal{LM}}$}}
\newcommand{\Q}{\mbox{${\mathcal{QM}}$}}
\newcommand{\RM}{\mbox{${\mathcal{RM}}$}}
\newcommand{\R}{\mbox{${\mathcal R}$}}
\newcommand{\Sy}{\mbox{${\mathcal S}$}}
\newcommand{\T}{\mbox{${\mathcal T}$}}
\newcommand{\U}{\mbox{${\mathcal U}_{loc}$}}
\newcommand{\UL}{\mbox{${\mathcal{UL}_{loc}}$}}
\newcommand{\UR}{\mbox{${\mathcal{UR}_{loc}}$}}
\newcommand{\ext}{\mbox{${\operatorname{ext}}$}}

\newcommand{\Ball}{\mbox{${\operatorname{Ball}}$}}

\newcommand{\TER}{\mbox{${\operatorname{TER}}$}}
\newcommand{\WTER}{\mbox{${\operatorname{WTER}}$}}
\begin{document}
\title[Extreme points of the unit ball of an operator space]{Extreme points of the unit ball of an operator space}
\author{Masayoshi Kaneda*}
\address{Department of Mathematics, University of California, Irvine, CA 92697-3875 U.S.A.}
\email{mkaneda@math.uci.edu, {\em URL}: http://www.math.uci.edu/$\thicksim$mkaneda/}
\address{Current Address: Department of Mathematics, The University of Mississippi, University, MS 38677 U.S.A.}
\date{\today}
\thanks{THIS PAPER IS A REVISION AND AN ENLARGEMENT OF THE AUTHOR'S MANUSCRIPT TITLED ``EXTREME POINTS OF THE UNIT BALL OF A QUASI-MULTIPLIER SPACE'' WHICH HAD BEEN CIRCULATED SINCE 2004}
\thanks{* The author was supported by a research fund from Department of Mathematics, University of California, Irvine}
\maketitle
\begin{center}Dedicated to the memory of Gert~K.~Pedersen\end{center}
\begin{abstract}We study extreme points of the unit ball of an operator space by introducing the new notion (approximate) ``quasi-identities''. Then we characterize an operator algebra with a contractive approximate quasi- (respectively, left, right, two-sided) identity in terms of quasi-multipliers and extreme points. Furthermore, we give a very neat necessary and sufficient condition for a given operator space to become a $C^*$-algebra or a one-sided ideal in a $C^*$-algebra in terms of quasi-multipliers. An extreme point is also used to show that any TRO with predual can be decomposed to the direct sum of a two-sided ideal, a left ideal, and a right ideal in some von Neumann algebra.
\end{abstract}
\section{Introduction.}\label{section: intro}
Quasi-multipliers of operator spaces were introduced by V.~I.~Paulsen (\cite{KP2}~Definition~2.2) in late 2002 as a natural variant of one-sided multipliers of operator spaces which had been introduced and studied by D.~P.~Blecher around 1999 (\cite{B4}). However, the significant correspondence of quasi-multipliers and operator algebra products was discovered by the author (\cite{KP2}~Theorem~2.6) in early 2003. That is, for a given operator space $X$, the possible operator algebra products that $X$ can be equipped with are precisely the bilinear mappings on $X$ that are implemented by contractive quasi-multipliers. Moreover, in \cite{Ka}~Theorem~3.3.1 and \cite{K}~Theorem~4.1, the author gave a geometric characterization of operator algebra products in terms of only matrix norms using the Haagerup tensor product. These results were presented in the Great Plains Operator Theory Symposium (GPOTS) held at the University of Illinois at Urbana-Champaign in 2003. After the author's talk, G.~K.~Pedersen asked the author the question ``How can the extreme points of the unit ball of a quasi-multiplier space be characterized?'' This question gave the author a further direction to study about quasi-multipliers. Through investigation, it turned out that what should be characterized is the ``interlocking'' between quasi-multipliers and extreme points of the unit ball of the operator space, and not extreme points of the unit ball of quasi-multiplier space alone. For this sake, the author introduces the new notion (approximate) ``quasi-identities''.

In Section~\ref{section: pre} we briefly review a construction of injective envelopes and triple envelopes of an operator space, and recall the definition of quasi-multipliers and their correspondence to operator algebra products. Furthermore, we define important classes of extreme points: local isometries, local co-isometries, and local unitaries, which actually becomes isometries, co-isometries, and unitaries, respectively, in certain cases with certain embeddings.

In Section~\ref{section:alternative} we give alternative definitions of one-sided and quasi-multipliers which are used to characterize operator algebras with approximate identities in Section~\ref{section: quasi-id}.

In Section~\ref{section: quasi-id} we introduce the new notion: (approximate) ``quasi-identities'' for normed algebras. We see that at least in the operator algebra case, contractive (approximate) quasi-identities are natural generalization of contractive (approximate) one-sided identities. Then we characterize an operator algebra with a contractive (approximate) quasi- (respectively, left, right, two-sided) identity in terms of its associated quasi-multiplier and extreme points of the unit ball of (the weak$^*$-closure of) the underlying operator space.

In Section~\ref{section: C*} we give an operator space characterization of $C^*$-algebras and their one-sided ideals in a very clear manner in terms of quasi-multipliers.

Section~\ref{section: dual} is devoted to showing that if an operator space has an operator space predual, then so is its quasi-multiplier space.

In Section~\ref{section:decompositions} it is shown that any TRO with predual can be decomposed to the direct sum of a two-sided ideal, a left ideal, and a right ideal in some von Neumann algebra using an extreme point of the TRO.

This paper is a revision and an enlargement of the author's manuscript titled ``Extreme points of the unit ball of an quasi-multiplier space'' which had been circulated since 2004. The author thanks David P. Blecher for pointing out a gap in the initial manuscript.

Portion of the work was carried out while the author was a post-doctoral researcher at the University of California, Irvine. The author is grateful to Bernard Russo for his invitation, financial support, and warm hospitality.

As stated before, this work was motivated by Pedersen's question posed to the author. However, he passed away within a
year of his asking the question. We regret that our answer to his question did not make it while he was alive.
Henceforth the author would like to dedicate this paper as a requiem to the memory of Gert~K.~Pedersen.
\section{Preliminaries.}\label{section: pre}
We begin by recalling the construction of an injective envelope of an operator space due to Z.-J.~Ruan (\cite{R}, \cite{R1}) and M.~Hamana (\cite{H}, \cite{H4}, also see \cite{H5}), independently. The reader unfamiliar with this subject is referred to \cite{P}~Chapter~15, \cite{ER2}~Chapter~6, \cite{BP2}, or \cite{B4}, for example.

Let $X\subset\B(\Hi)$ be an operator space, and consider Paulsen's operator system$$\Sy_X:=\left[\begin{matrix}\C1_{\Hi}&X\\X^*&\C1_{\Hi}\\\end{matrix}\right]\subset\M_2(\B(\Hi)).$$One then takes a minimal (with respect to a certain ordering) completely positive $\Sy_X$-projection\footnote{An $\Sy_X$-projection is an idempotent that fixes each element in $\Sy_X$.} $\Phi$ on $\M_2(\B(\Hi))$, whose image Im$\Phi$ turns out to be an injective envelope $I(\Sy_X)$ of $\Sy_X$. By a well-known result of M.-D.~Choi and E.~G.~Effros (\cite{CE}), Im$\Phi$ is a unital $C^*$-algebra with the product $\odot$ (which is called the Choi-Effros product) defined by $\xi\odot\eta:=\Phi(\xi\eta)$ for $\xi,\eta\in$Im$\Phi$ and with other algebraic operations and norm taken to be the original ones in $\M_2(\B(\Hi))$. The $C^*$-algebra structure of $I(\Sy_X)$ does not depend on the particular embedding $X\subset\B(\Hi)$. By a well-known trick one may decompose to$$\Phi=\left[\begin{matrix}\psi_1&\phi\\
\phi^*&\psi_2\end{matrix}\right].$$Accordingly, one may write
\begin{equation}\label{eq:injenv}\text{Im}\Phi=I(\Sy_X)=\left[\begin{matrix}I_{11}(X)&I(X)\\I(X)^*&I_{22}(X)\end{matrix}
\right]\subset\M_2(\B(\Hi)),
\end{equation}where $I(X)$ is an injective envelope of $X$, and $I_{11}(X)$ and $I_{22}(X)$ are injective $C^*$-algebras (hence unital (See \cite{BP2}~Proposition~2.8.)). We denote the identities of $I_{11}(X)$ and $I_{22}(X)$ by $1_{11}$ and $1_{22}$, respectively. Note that the last inclusion in Expression~(\ref{eq:injenv}) is not as a subalgebra since the multiplication in $I(\Sy_X)$ and the multiplication in $\M_2(\B(\Hi))$ are not same in general. The new product $\odot$ induces a new product $\bullet$ between elements of $I_{11}(X)$, $I_{22}(X)$, $I(X)$, and $I(X)^*$. For instance, $x\bullet y^*=\psi_1(xy^*)$ for $x,y\in I(X)$. Note that the associativity of $\bullet$ is guaranteed by that of $\odot$.

The following property is often useful.
\begin{lemma}\label{lm: BP}{\em(Blecher-Paulsen \cite{BP2}~Corollary~1.3)}
\begin{enumerate}
\item If $a\in I_{11}(X)$, and if $a\bullet x=0,\;\forall x\in X$, then $a=0$.
\item If $b\in I_{22}(X)$, and if $x\bullet b=0,\;\forall x\in X$, then $b=0$.
\end{enumerate}
\end{lemma}One may write the $C^*$-subalgebra $C^*(\partial X)$ of Im$\Phi$ (with the new product) generated by$$\left[\begin{matrix}O&X\\O&O\end{matrix}\right]$$as$$C^*(\partial X)=\left[\begin{matrix}\E(X)&\T(X)\\\T(X)^*&\F(X)\end{matrix}\right]\subset\left[\begin{matrix}I_{11}(X)&I(X)\\I(X)^*&
I_{22}(X)\end{matrix}\right],$$where $\T(X)$ is a triple envelope of $X$, i.e., a ``minimum'' TRO that contains $X$ completely isometrically. Here an operator space $X$ being a {\bf ternary ring of operators} ({\bf TRO} for short) or a {\bf triple system} means that there is a complete isometry $\iota$ from $X$ into a $C^*$-algebra such that $\iota(x)\iota(y)^*\iota(z)\in\iota(X),\;\forall x,y,z\in X$.

We call the embedding $i:X\to\left[\begin{matrix}O&X\\O&O\end{matrix}\right]\subset C^*(\partial X)\subset I(\Sy_X)$ the {\bf \v{S}ilov embedding} of $X$, and often denote $\left[\begin{matrix}0&x\\0&0\end{matrix}\right]\;(x\in X)$ simply by $x$. Similarly, we often write $X$ for $\left[\begin{matrix}O&X\\O&O\end{matrix}\right]$, and $1_{11}$ for $\left[\begin{matrix}1_{11}&0\\0&0\end{matrix}\right]$, etc. The involution in $\B(\Hi)$ induces an involution in
$\M_2(\B(\Hi))$ in an obvious way, and we still denote by $*$. For example, for $x\in X$, $\left[\begin{matrix}0&x\\0&0\end{matrix}\right]^*=\left[\begin{matrix}0&0\\x^*&0\end{matrix}\right]$.

In this paper, all operator spaces are assumed to be norm-closed. Whenever an infinite-dimensional vector space is involved in a product, we take the norm closure of the linear span. For instance, $X\bullet z\bullet X:=\overline{\text{span}}\{x\bullet z\bullet y;\;x,y\in X\}$, where $z\in I(X)^*$.

Now we are ready to recall the definition of one-sided and quasi-multipliers. We remark that the one-sided multipliers were first introduced by D.~P.~Blecher in \cite{B4}. The following definition (Items~(1)~and~(2)) is an equivalent but more manageable version in \cite{BP2}.
\begin{definition}\label{def: qm}Let $X$ be an operator space.
\begin{enumerate}
\item{\em (Blecher-Paulsen \cite{BP2}~Definition~1.4)} The {\bf left multiplier algebra\footnote{In \cite{BP2} it is denoted by $IM_l(X)$.}} of $X$ is the operator algebra$$\LM(X):=\{a\in I_{11}(X)\;;\;a\bullet X\subset X\}.$$We call an element of $\LM(X)$ a {\bf left multiplier} of $X$.
\item The {\bf right multiplier algebra\footnote{In \cite{BP2} it is denoted by $IM_r(X)$.}} of $X$ is the operator algebra$$\RM(X):=\{b\in I_{22}(X)\;;\;X\bullet b\subset X\}.$$We call an element of $\RM(X)$ a {\bf right multiplier} of $X$.
\item{\em (Paulsen \cite{KP2}~Definition~2.2)} The {\bf quasi-multiplier space} of $X$ is the operator space$$\Q(X):=\{z\in I(X)^*\;;\;X\bullet z\bullet X\subset X\}.$$We call an element of $\Q(X)$ a {\bf quasi-multiplier} of $X$.
\end{enumerate}
\end{definition}The following theorem characterizes operator algebra products in terms of quasi-multipliers and matrix norms. Especially, (iii) tells us that the operator algebra products (algebraic property) a given operator space can be equipped with are completely determined only by its underlying matrix norm structure (geometric property), which can be regarded as the ``quasi'' version of the $\tau$-trick theorem by Blecher-Effros-Zarikian (\cite{BEZ}~Theorem~1.1, Theorem~4.6).
\begin{theorem}\label{th: kaneda}{\em (Kaneda \cite{Ka}~Theorem~3.3.1, \cite{K}~Theorem~4.1)} Let $X$ be a non-zero operator space with a bilinear mapping $\varphi: X\times X\to X$, and let $I(\Sy_X)$ be as above and 1 be its identity. We regard $X$ as a subspace of $I(\Sy_X)$ by the \v{S}ilov embedding as explained above. Let$$\begin{matrix}&\M_2(I(\Sy_X)\stackrel{h}{\otimes}I(\Sy_X))&&\M_2(X)\\&\cup&&\cup\\\Gamma_{\varphi}:&\left[\begin
{matrix}X\stackrel{h}{\otimes}\C1&X\stackrel{h}{\otimes}X\\O&\C1\stackrel{h}{\otimes}X\end{matrix}\right]&\to&
\left[\begin{matrix}X&X\\O&X\end{matrix}\right]\end{matrix}$$be defined by$$\Gamma_{\varphi}\left(\left[\begin{matrix}x_1\otimes1&x\otimes y\\0&1\otimes x_2\end{matrix}\right]\right):=\left[\begin{matrix}x_1&\varphi(x,y)\\0&x_2\end{matrix}\right]$$and their linear extension and norm closure, where $\stackrel{h}{\otimes}$ is the Haagerup tensor product. Then, the following are equivalent:
\begin{enumerate}
\item[(i)]$(X,\varphi)$ is an abstract operator algebra (i.e., there is a completely isometric homomorphism from $X$ into a concrete operator algebra, hence, in particular, $\varphi$ is associative);
\item[(ii)]There exists a $z\in\Ball(\Q(X))$ \footnote{For a normed space $Y$, $\Ball(Y):=\{y\in Y\;;\;\|y\|\le1\}$.} such that $\forall x,y\in X$, $\varphi(x,y)=x\bullet z\bullet y$;
\item[(iii)]$\Gamma_{\varphi}$ is completely contractive.
\end{enumerate}Moreover, such a $z$ is unique.

When these conditions hold, we denote $\varphi$ by $m_z$, and call $(X,m_z)$ the operator algebra {\bf corresponding to} the quasi-multiplier $z$, or the {\bf algebrization} of $X$ by the quasi-multiplier $z$. On the other hand, for a given operator algebra $\A$, we call $z$ given in (ii) the quasi-multiplier {\bf associated with} $\A$.
\end{theorem}
We denote the set of the extreme points of the unit ball of the quasi-multiplier space of an operator space $X$ by $\ext(\Ball(\Q(X)))$. The following are particularly important subsets of $\ext(\Ball(\Q(X)))$. We will see in Corollary~\ref{cor: subext} that these are actually subsets of $\ext(\Ball(\Q(X)))$.
\begin{definition}\label{def:local}Let $X$ be an operator space, and let $S$ be a subset of $I(X)^*$.
\begin{enumerate}
\item $\UL(S):=\{z\in S\;;\;z^*\bullet z=1_{11}\}$.
\item $\UR(S):=\{z\in S\;;\;z\bullet z^*=1_{22}\}$.
\item $\U(S):=\UL(S)\cap\UR(S)$.
\end{enumerate}We call an element of $\UL(I(X)^*)$ (respectively, $\UR(I(X)^*)$, $\U(I(X)^*)$) a {\bf local isometry} (or, {\bf local left unitary}) (respectively, {\bf local co-isometry} (or, {\bf local right unitary}), {\bf local\footnote{Matthew~Neal called them ``local'' when the author presented this manuscript to him.} unitary}).
\end{definition}Item (3) of the following proposition tells us that with certain embeddings, local isometries (respectively, local co-isometries, local unitaries) actually become isometries (respectively, co-isometries, unitaries).
\begin{lemma}\label{lm:embed}Let $X$ be an operator space.
\begin{enumerate}
\item If $\UL(I(X)^*)\ne\varnothing$, then there exits a commutative diagram$$\begin{array}{ccc}I_{11}(X)&\stackrel{\sigma_1}{\longrightarrow}&I(X)\\\sigma_1^*\downarrow&
\circlearrowright&\downarrow\rho_1\\I(X)^*&\stackrel{\rho_1^*}{\longrightarrow}&I_{22}(X)\end{array}$$such that $\rho_1^*(x^*):=\rho_1(x)^*,\;\forall x\in I(X)$, $\sigma_1^*(a):=\sigma_1(a^*)^*,\;\forall a\in I_{11}(X)$, and $\rho_1$, $\sigma_1$ (hence $\rho_1^*$, $\sigma_1^*$) are complete isometries, and $\rho_1\circ\sigma_1$ (hence $\rho_1^*\circ\sigma_1^*$) is a $*$-monomorphism, and $\forall z\in\UL(I(X)^*)$, $\rho_1^*(z)$ is a partial isometry in the $C^*$-algebra $I_{22}(X)$.
\item If $\UR(I(X)^*)\ne\varnothing$, then there exits a commutative diagram$$\begin{array}{ccc}I_{11}(X)&\stackrel{\sigma_2}{\longleftarrow}&I(X)\\\sigma_2^*\uparrow&\circlearrowleft& \uparrow\rho_2\\I(X)^*&\stackrel{\rho_2^*}{\longleftarrow}&I_{22}(X)\end{array}$$such that $\rho_2^*(b):=\rho_2(b^*)^*,\;\forall b\in I_{22}(X)$, $\sigma_2^*(x^*):=\sigma_2(x)^*,\;\forall x\in I(X)$, and $\rho_2$, $\sigma_2$ (hence $\rho_2^*$, $\sigma_2^*$) are complete isometries, and $\sigma_2\circ\rho_2$ (hence $\sigma_2^*\circ\rho_2^*$) is a $*$-monomorphism, and $\forall z\in\UR(I(X)^*)$, $\sigma_2^*(z)$ is a partial isometry in the $C^*$-algebra $I_{11}(X)$.
\item If $\U(I(X)^*)\ne\varnothing$, then in (1) and (2), one can take $\rho_1$, $\sigma_1$, $\rho_2$, $\sigma_2$ (hence $\rho_1^*$, $\sigma_1^*$, $\rho_2^*$, $\sigma_2^*$) to be onto such that $\rho_2=\rho_1^{-1}$, $\sigma_2=\sigma_1^{-1}$ (hence $\rho_2^*=(\rho_1^*)^{-1}$, $\sigma_2^*=(\sigma_1^*)^{-1}$). Moreover,
    \begin{enumerate}
    \item $\forall z\in\UL(I(X)^*)$, $\sigma_2^*(z)=(\sigma_1^*)^{-1}(z)$ and $\rho_1^*(z)=(\rho_2^*)^{-1}(z)$ are isometries in the $C^*$-algebras $I_{11}(X)$ and $I_{22}(X)$, respectively.
    \item $\forall z\in\UR(I(X)^*)$, $\sigma_2^*(z)=(\sigma_1^*)^{-1}(z)$ and $\rho_1^*(z)=(\rho_2^*)^{-1}(z)$ are co-isometries in the $C^*$-algebras $I_{11}(X)$ and $I_{22}(X)$, respectively.
    \item $\forall z\in\U(I(X)^*)$, $\sigma_2^*(z)=(\sigma_1^*)^{-1}(z)$ and $\rho_1^*(z)=(\rho_2^*)^{-1}(z)$ are unitaries in the $C^*$-algebras $I_{11}(X)$ and $I_{22}(X)$, respectively.
    \end{enumerate}
\end{enumerate}
\end{lemma}
\begin{proof}Once we define mappings $\rho_1$, $\rho_2$, $\sigma_1$, and $\sigma_2$ as follows, then the assertions are straightforward.

$\underline{\text{(1):}}$ Pick $z_1\in\UL(I(X)^*)$, and define $\rho_1(x):=z_1\bullet x,\forall x\in I(X)$; $\sigma_1(a):=a\bullet z_1^*,\forall a\in I_{11}(X)$.

$\underline{\text{(2):}}$ Pick $z_2\in\UR(I(X)^*)$, and define $\rho_2(b):=z_2^*\bullet b,\forall b\in I_{22}(X)$; $\sigma_2(x):=x\bullet z_2,\forall x\in I(X)$.

$\underline{\text{(3):}}$ Pick $z_3\in\U(I(X)^*)$, and define $\rho_1$, $\rho_2$, $\sigma_1$, and $\sigma_2$, as (1) and (2).
\end{proof}
\begin{lemma}\label{lm:BK}Let $\A$ be a nonzero operator algebra.
\begin{enumerate}
\item If $\A$ has a contractive approximate left identity, then $\UL(\Q(\A))\ne\varnothing$.
\item If $\A$ has a contractive approximate right identity, then $\UR(\Q(\A))\ne\varnothing$.
\item If $\A$ has a contractive approximate two-sided identity, then $\U(\Q(\A))\ne\varnothing$.
\end{enumerate}
\end{lemma}
\begin{proof}(2) follows from noticing that in the proof of Theorem~2.3 in \cite{BK}, it is seen that $v^*$ is the quasi-multiplier associated with $\A$ and  $v^*\in\UR(\Q(\A))$, although the terminology ``quasi-multiplier'' does not appear there. (1) is similar by symmetry. (3) follows from (1) and (2) together with the uniqueness of a quasi-multiplier associated with $\A$ (Theorem~\ref{th: kaneda} of the present paper).
\end{proof}It follows from Lemma~\ref{lm:BK} together with Lemma~\ref{lm:embed} that if $\A$ is an operator algebra with a contractive approximate right identity, then $I(\A)$ and $I(\A)^*$ are embedded in $I_{11}(\A)$ completely isometrically, and $I_{22}(\A)$ is embedded in $I_{11}(\A)$ $*$-monomorphically. A part of this fact is already seen in \cite{BK}~Theorem~2.3. Similar embeddings hold when $\A$ has a contractive approximate left identity. Also note that in \cite{Ka}~Lemma~3.2.2 and \cite{K}~Lemma~3.4 assuming that $\A$ has a contractive approximate two-sided identity we embedded $I(\A)$ and $I(\A)^*$ in $I_{11}(\A)$, accordingly we showed that our definition of quasi-multipliers (Definition~\ref{def: qm}) coincides with the classical ones for $C^*$-algebras (\cite{Pe}~Section 3.12) in the sense that they are completely isometrically quasi-isomorphic (\cite{Ka}~Theorem 3.2.3 and \cite{K}~Theorem~3.5).

Recall that if $\A$ is an operator algebra with a contractive approximate two-sided identity, then its injective envelope $I(\A)$ is a unital $C^*$-algebra which contains $\A$ as a subalgebra (See \cite{BL}~Corollary~4.2.8~(1) for example.).
\begin{definition}\label{def:unital}
\begin{enumerate}
\item $\LM^1(\A):=\{a\in I(\A)\;;\;a\A\subset\A\}$.
\item $\RM^1(\A):=\{b\in I(\A)\;;\;\A b\subset\A\}$.
\item $\Q^1(\A):=\{z\in I(\A)\;;\;\A z\A\subset\A\}$.
\end{enumerate}
\end{definition}$\LM^1(X)$, $\RM^1(X)$, and $\Q^1(X)$ are equivalent to $\LM(X)$, $\RM(X)$, and $\Q(X)$, respectively in the sense of the following lemma.
\begin{lemma}\label{lm:unital}Let $\A$ be a nonzero operator algebra with a contractive approximate two-sided identity. Then the following assertions hold.
\begin{enumerate}
\item There is a multiplicative complete isometry $\lambda$ from $\LM^1(\A)$ onto $\LM(\A)$ such that $a x=\lambda(a)\bullet x,\;\forall a\in\LM^1(\A),\forall x\in\A$.
\item There is a multiplicative complete isometry $\rho$ from $\RM^1(\A)$ onto $\RM(\A)$ such that $x b=x\bullet\rho(b),\;\forall b\in\RM^1(\A),\forall x\in\A$.
\item There is a complete isometry $\kappa$ from $\Q^1(\A)$ onto $\Q(\A)$ such that $xzy=x\bullet\kappa(z)\bullet y,\;\forall z\in\Q^1(\A),\forall x,y\in\A$.
\end{enumerate}
\end{lemma}
\begin{proof}To see (3), first note that by Lemma~\ref{lm:BK}~(3) and Lemma~\ref{lm:embed}~(3), $I_{11}(\A)$ is an injective envelope of $\A$. By the uniqueness of an injective envelope $C^*$-algebra up to $*$-isomorphism that fixes each element of $\A$, $I(\A)$ in the definition of $\Q^1(\A)$ can be taken to be $I_{11}(\A)$. Now the assertion follows from \cite{Ka}~Lemma~3.2.2 (or \cite{K}~Lemma~3.4). Items (1) and (2) are similar by developing a lemma corresponding to \cite{Ka}~Lemma~3.2.2 (or \cite{K}~Lemma~3.4). The details are left to the reader.
\end{proof}We close this preliminary section with Kadison's characterization of the extreme points of the unit ball of a $C^*$-algebra (\cite{Kad}~Theorem~1). We use the following version in Pedersen's book (\cite{Pe}~Proposition~1.4.8) or Sakai's book
(\cite{S2}~Proposition~1.6.5). This motivated our definition of quasi-identities and plays a key role in the proof of the characterization theorems (Theorems~\ref{th:extinj}~and~\ref{th:main}) and ideal decompositions (Theorem~\ref{th: ideal decomposition}).
\begin{lemma}\label{lm: kadison}{\em (Kadison)} Let $\A$ be a $C^*$-algebra, and let $p$, $q$ be orthogonal projections in $\A$. Then the following are equivalent:
\begin{enumerate}
\item[(i)]$x\in p\A q$ is an extreme point of $\Ball(p\A q)$;
\item[(ii)]$(\tilde{p}-xx^*)p\A q(\tilde{q}-x^*x)=\{0\}$ for some orthogonal projections $\tilde{p}$ and $\tilde{q}$ in $\A$ such that $\tilde{p}\ge p$ and $\tilde{q}\ge q$;
\item[(iii)]$(\tilde{p}-xx^*)p\A q(\tilde{q}-x^*x)=\{0\}$ for all orthogonal projections $\tilde{p}$ and $\tilde{q}$ in $\A$ such that $\tilde{p}\ge p$ and $\tilde{q}\ge q$.
\end{enumerate}
In this case, $x$ is a partial isometry in $\A$.
\end{lemma}The following corollary is immediate from the lemma above.
\begin{corollary}\label{cor: subext}Let $X$ be an operator space. Then $\UL(\Q(X))$, $\UR(\Q(X))$, and $\U(\Q(X))$ are subsets of $\ext(\Ball(I(X)^*))$, hence subsets of $\ext(\Ball(\Q(X)))$.
\end{corollary}
\section{Alternative definitions of one-sided and quasi-multipliers}\label{section:alternative}
In this section we give alternative definitions of one-sided and quasi-multipliers of an operator space $X$ which are equivalent to the ones presented in Definition \ref{def: qm}.

First, we give a definition of one-sided and quasi-multipliers of $X$ using the second dual of $C^*(\partial X)$. Denote\footnote{To avoid confusion with the adjoint and also to distinguish $\Q''(\A)$ from $\Q^{**}(\A)$ in Item~(I) on page~86 of \cite{Ka} and Item~(I) on page~351 of \cite{K}, in this paper we denote the second dual by the double primes instead of the double stars.} the second dual of $\E(X)$, $\F(X)$, and $\T(X)$ by $\E(X)''$, $\F(X)''$, and $\T(X)''$, respectively, and we regard them as the corners of the second dual of $C^*(\partial X)$ in the usual way:$$C^*(\partial X)''=\left[\begin{matrix}\E(X)''&\T(X)''\\{\T(X)^*}''&\F(X)''\end{matrix}\right].$$The Arens product on $C^*(\partial X)''$ induces a product between elements of $\E(X)''$, $\F(X)''$, $\T(X)''$, and ${\T(X)^*}''$, which is an extension of $\bullet$ defined in Section~\ref{section: pre} and is still denoted by $\bullet$. Denote by $1_{\E}$ and $1_{\F}$ the identity of the $W^*$-algebras $\E(X)''$ and $\F(X)''$, respectively. Let $\widehat{\quad}:C^*(\partial X)\to C^*(\partial X)''$ be the canonical embedding.
\begin{definition}\label{def:alternative}
\begin{enumerate}
\item$\LM''(X):=\{a\in\E(X)''\;;\;a\bullet\widehat{X}\subset\widehat{X}\}$.
\item$\RM''(X):=\{b\in\F(X)''\;;\;\widehat{X}\bullet b\subset\widehat{X}\}$.
\item$\Q''(X):=\{z\in{\T(X)^*}''\;;\;\widehat{X}\bullet z\bullet\widehat{X}\subset\widehat{X}\}$.
\end{enumerate}
\end{definition}
$\LM''(X)$, $\RM''(X)$, and $\Q''(X)$ are equivalent to $\LM(X)$, $\RM(X)$, and $\Q(X)$, respectively in the sense of
the following proposition which we will prove shortly.
\begin{proposition}\label{pr:alternative}
\begin{enumerate}
\item There is a multiplicative completely isometry $\lambda_1$ from $\LM(X)$ onto $\LM''(X)$ such that $\widehat{a\bullet x}=\lambda_1(a)\bullet\hat{x},\;\forall a\in\LM(X),\forall x\in X$.
\item There is a multiplicative completely isometry $\rho_1$ from $\RM(X)$ onto $\RM''(X)$ such that $\widehat{x\bullet b}=\hat{x}\bullet\rho_1(b),\;\forall b\in\RM(X),\forall x\in X$.
\item There is a completely isometry $\kappa_1$ from $\Q(X)$ onto $\Q''(X)$ such that $\widehat{x\bullet z\bullet y}=\hat{x}\bullet\kappa_1(z)\bullet\hat{y},\;\forall z\in\Q(X),\forall x,y\in X$.
\end{enumerate}
\end{proposition}Next, we give a definition of one-sided and quasi-multipliers of $X$ using a representation of $C^*(\partial X)$ on a Hilbert space. Represent $C^*(\partial X)$ by a $*$-monomorphism $\pi$ on the direct sum of Hilbert spaces $\Hi_1$ and $\Hi_2$ nondegenerately so that $[\T(X)\Hi_2]=\Hi_1$ and $[\T(X)^*\Hi_1]=\Hi_2$.\footnote{$[\T(X)\Hi_2]:=\overline{\operatorname{span}}\{x\xi\;;\;x\in\T(X),\xi\in\Hi_2\}$.} Denote by $1_{\Hi_1}$ and $1_{\Hi_2}$ the orthogonal projections onto $\Hi_1$ and $\Hi_2$, respectively.
\begin{definition}\label{def:alternative2}
\begin{enumerate}
\item$\LM_{\pi}(X):=\{a\in\B(\Hi_1)\;;\;a\pi(X)\subset\pi(X)\}$.
\item$\RM_{\pi}(X):=\{b\in\B(\Hi_2)\;;\;\pi(X)b\subset\pi(X)\}$.
\item$\Q_{\pi}(X):=\{z\in\B(\Hi_1,\Hi_2)\;;\;\pi(X)z\pi(X)\subset\pi(X)\}$.
\end{enumerate}
\end{definition}
$\LM_{\pi}(X)$, $\RM_{\pi}(X)$, and $\Q_{\pi}(X)$ are equivalent to $\LM(X)$, $\RM(X)$, and $\Q(X)$, respectively in the sense of the following proposition.
\begin{proposition}\label{pr:alternative2}
\begin{enumerate}
\item There is a multiplicative completely isometry $\lambda_2$ from $\LM(X)$ onto $\LM_{\pi}(X)$ such that $\pi(a\bullet x)=\lambda_2(a)\pi(x),\;\forall a\in\LM(X),\forall x\in X$.
\item There is a multiplicative completely isometry $\rho_2$ from $\RM(X)$ onto $\RM_{\pi}(X)$ such that $\pi(x\bullet b)=\pi(x)\rho_2(b),\;\forall b\in\RM(X),\forall x\in X$.
\item There is a completely isometry $\kappa_2$ from $\Q(X)$ onto $\Q_{\pi}(X)$ such that $\pi(x\bullet z\bullet y)=\pi(x)\kappa_2(z)\pi(y),\;\forall z\in\Q(X),\forall x,y\in X$.
\end{enumerate}
\end{proposition}
{\em Proof of Propositions~\ref{pr:alternative}~and~\ref{pr:alternative2}.} To see Proposition~\ref{pr:alternative2}~(3), first note that $I(\Sy_X)$ is also an injective envelope $C^*$-algebra of
$C^*(\partial X)$. Since a $C^*$-algebra has a contractive approximate two-sided identity, there is a complete isometry
$\kappa$ from $\Q^1(C^*(\partial X))$ onto $\Q^{\pi}(C^*(\partial X))$ \footnote{This $\Q^{\pi}(C^*(\partial X))$ with ``superscript $\pi$'' is as defined in Item~(II) on page~86 of \cite{Ka} and Item~(II) on page~351 of \cite{K}, and is different from $\Q_{\pi}(C^*(\partial X))$ with ``subscript $\pi$'' defined in Definition~\ref{def:alternative2} of the present paper, although they are quasi-isomorphic in the sense of \cite{Ka}~Definition~3.1.1~(2) and \cite{K}~Definition~2.1~(2).} such that $\pi(\xi\zeta\eta)=\pi(\xi)\kappa(\zeta)\pi(\eta),\forall\zeta\in\Q^1(C^*(\partial X)),\forall\xi,\eta\in C^*(\partial X)$ by Lemma~\ref{lm:unital}~(3) together with \cite{Ka}~Theorem~3.2.3 (or \cite{K}~Theorem~3.5). The restriction of $\kappa$ to $\Q(X)$ gives $\kappa_2$. Proposition~\ref{pr:alternative}~(3) is similar. (1) and (2) of Propositions~\ref{pr:alternative}~and~\ref{pr:alternative2} are also similar, but use \cite{B}~Theorem~6.1 instead of \cite{Ka}~Theorem~3.2.3 (or \cite{K}~Theorem~3.5), and the equivalence of $\LM(C^*(\partial X))$ and $M_l(C^*(\partial X))$ (\cite{BP2}~Theorem~1.9~(i)).\begin{flushright}$\square$\end{flushright}

Note that it is possible to write Theorem~\ref{th: kaneda} using these alternative definitions.

Finally, we define the following sets.
\begin{definition}\label{def:local-alternative}Let $X$ be an operator space.
\begin{enumerate}
\item Let $S$ be a subset of ${\T(X)^*}''$.
\begin{enumerate}
\item$\UL(S):=\{z\in S\;;\;z^*\bullet z=1_{\E}\}$.
\item$\UR(S):=\{z\in S\;;\;z\bullet z^*=1_{\F}\}$.
\item$\U(S):=\UL(S)\cap\UR(S)$.
\end{enumerate}
\item Let $S$ be a subset of $\B(\Hi_1,\Hi_2)$.
\begin{enumerate}
\item$\UL(S):=\{z\in S\;;\;z^*z=1_{\Hi_1}\}$.
\item$\UR(S):=\{z\in S\;;\;zz^*=1_{\Hi_2}\}$.
\item$\U(S):=\UL(S)\cap\UR(S)$.
\end{enumerate}
\end{enumerate}
\end{definition}One may rewrite Lemma~\ref{lm:embed} and Lemma~\ref{lm:BK} using Definition~\ref{def:local-alternative}. The details are left to the reader.

The following corollary immediately follows from Kadison's Theorem (Lemma~\ref{lm: kadison}).
\begin{corollary}\label{cor:subext}
\begin{enumerate}
\item$\UL(\Q''(X))$, $\UR(\Q''(X))$, and $\U(\Q''(X))$ are subsets of $\ext(\Ball(\T(X)^*))$, hence subsets of $\ext(\Ball(\Q''(X))$.
\item$\UL(\Q_{\pi}(X))$, $\UR(\Q_{\pi}(X))$, and $\U(\Q_{\pi}(X))$ are subsets of $\ext(\Ball(\B(\Hi,\K)))$, hence subsets of $\ext(\Ball(\Q_{\pi}(X))$.
\end{enumerate}
\end{corollary}Hereafter, we omit the symbol $\iota$ or $\pi$, and we regard $C^*(\partial X)$ as a $C^*$-subalgebra of $C^*(\partial X)''$ or $\B(\Hi)$. Also we omit the symbol $\odot$ or $\bullet$ unless there is a possibility of
confusion.
\section{Quasi-identities and characterization theorems}\label{section: quasi-id}
Throughout this section, the following elementary lemma which follows from the polarization identity is useful.
\begin{lemma}\label{lm: idempotent}
\begin{enumerate}
\item Let $a\in\B(\Hi)$. If $a^2=a$ and $\|a\|\le1$, then $a^*=a$, i.e., $a$ is an orthogonal projection.
\item Let $p\in\B(\Hi)$ be an orthogonal projection, i.e., $p=p^*=p^2$, and let $b,c\in\B(\Hi)$ such that $c^*b=p$ and $\|b+c\|\le2$. Then $\ker p\subset\ker b\cap\ker c$ if and only if $b=c$. In this case, $\ker p=\ker b=\ker c$
\end{enumerate}
\end{lemma}
\begin{proof}$\underline{\text{(1):}}$ Let $\xi\in\Hi$. Then by the polarization identity, $\|a\xi\|^2=<a\xi,a\xi>=<a\xi,a^2\xi>=<a^*a\xi,a\xi>=\frac{1}{4}(\|(a^*a+a)\xi\|^2-\|(a^*a-a)\xi\|^2)\le\|a\xi\|^
2-\frac{1}{4}\|(a^*a-a)\xi\|^2$. Since $\xi\in\Hi$ is arbitrary, $a=a^*a=a^*$.

$\underline{\text{(2):}}$ Assume that $\ker p\subset\ker b\cap\ker c$, and let $\eta\in\Hi$. Then by the polarization identity, $\|p\eta\|^2=<c^*bp\eta,p\eta>=<bp\eta,cp\eta>=\frac{1}{4}(\|(b+c)p\eta\|^2-\|(b-c)p\eta\|^2)\le\|p\eta\|^2-\frac{1}{4}
\|(b-c)p\eta\|^2$ from which it follows that $b=c$. The converse direction and the last assertion are obvious and very basic facts.
\end{proof}
We introduce the new notion ``(approximate) quasi-identities''.
\begin{definition}\label{def: qi}
\begin{enumerate}
\item Let $\R$ be a ring. A {\bf quasi-identity} of $\R$ is an element $e\in\R$ such that$$r=er+re-ere,\quad\forall r\in\R.$$
\item Let $\A$ be a normed algebra. An {\bf approximate quasi-identity} of $\A$ is a net $\{e_\alpha\}\subset\A$ such that$$a=\lim_{\alpha\to\infty}(e_{\alpha}a+ae_{\alpha}-e_{\alpha}ae_{\alpha}),\quad\forall a\in\A.$$
\end{enumerate}
\end{definition}It is quite essential in the definition of an approximate quasi-identity that the limit is taken ``at once''. In fact, a bounded approximate left identity $\{e_{\alpha}\}$ of a normed algebra $\A$ is easily seen to be an approximate quasi-identity. However, $\lim_{\alpha}ae_{\alpha}$ need not exist for all $a\in\A$ as is seen in Example~\ref{ex:B}. Despite this fact, we restrict our characterization of contractive approximate quasi-identities to the case that both $\lim_{\alpha}e_{\alpha}a$ and $\lim_{\alpha}ae_{\alpha}$ exist for all $a\in\A$ (See the proof of Theorem~\ref{th:main}~(1).). This somewhat unpleasant point is due to the fact that an operator algebra product is weak$^*$-continuous with respect to each factor ``separately''.
\begin{proposition}\label{pr:qi}
\begin{enumerate}
\item A separable normed algebra with an approximate quasi-identity admits an approximate quasi-identity which is a sequence.
\item A finite-dimensional normed algebra with a bounded approximate quasi-identity actually has a quasi-identity.
\end{enumerate}
\end{proposition}
\begin{proof}Item~(1) can be proved in a similar way to showing a separable $C^*$-algebra admits an approximate identity which is a sequence (See \cite{M}~Remark~3.1.1 for example.), and the details are left to the reader. To see (2), let $\{e_{\alpha}\}$ be a bounded approximate quasi-identity of a finite-dimensional normed algebra $\A$, and let $e$ be an accumulation point of $\{e_{\alpha}\}$ in $\A$. Then one can take a subnet $\{e_{\alpha_n}\}$ (which can be a sequence) such that $\lim_ne_{\alpha_n}=e$. Therefore $\forall a\in\A$, $\|a-(ea+ae-eae)\|\le\|a-(e_{\alpha_n}a+ae_{\alpha_n}-e_{\alpha_n}ae_{\alpha_n})\|+\|e_{\alpha_n}-e\|\|a\|
+\|a\|\|e_{\alpha_n}-e\|+\|e_{\alpha_n}-e\|\|a\|\|e_{\alpha_n}\|+\|e\|\|a\|\|e_{\alpha_n}-e\|\to0$ as $n\to\infty$.
\end{proof}Identities, left identities, right identities of rings are quasi-identities. We will see in Proposition~\ref{pr:unique} that in the operator algebra case, a contractive quasi-identity is unique if it exists. Moreover, it is necessarily idempotent and Hermitian (if the operator algebra is embedded in a $C^*$-algebra by a multiplicative complete isometry).

Bounded approximate left (respectively, right, two-sided) identities of normed algebras are approximate quasi-identities. Many normed algebras do not have an (approximate) two-sided or one-sided identity, but do have an (approximate) quasi-identity. Perhaps the following illustrates a typical situation. Let $\A$ be a normed algebra which has a bounded left approximate identity $\{e_\alpha\}$ but does not have a right approximate identity, and let $\Br$ be a normed algebra which has a bounded right approximate identity $\{f_\beta\}$ but does not have a left approximate identity. Then $\A\stackrel{p}{\oplus}\Br$ with $1\le p\le\infty$ has neither left nor right approximate identity, but does have a bounded approximate quasi-identity $\{e_\alpha\oplus f_\beta\}_{(\alpha,\beta)}$, where $\{(\alpha,\beta)\}$ is a directed set by the ordering defined by ``$(\alpha_1,\beta_1)\le(\alpha_2,\beta_2)$ if and only if $\alpha_1\le\alpha_2$ and $\beta_1\le\beta_2$''. But if we can always decompose a normed algebra to the direct sum of two normed algebras one of which has a left approximate identity and the other has a right approximate identity, then it is not so meaningful to define (approximate) quasi-identities since we can always reduce to the case of normed algebras with a one-sided approximate identity. We thank Takeshi~Katsura for asking for such an example that cannot be decomposed to the direct sum of two normed algebras with a one-sided identity. Here is an example: Let $\A$ be the subalgebra of $\M_3(\C)$ supported on the $(1,1)$-, $(1,2)$-, $(1,3)$-, $(2,3)$-, and $(3,3)$-entries only, where $\M_3(\C)$ is equipped with the usual matrix operations. Then $\A$ has neither left nor right identity, and cannot be decomposed to the direct sum of any two algebras, but does have a quasi-identity $E_1+E_3$, where $E_i$ denotes the matrix whose $(i,i)$-entry is $1$ and all other entries are $0$'s.

The following proposition convinces us that the notion of (approximate) quasi-identities is natural and in a certain sense ``minimal'' generalization of (approximate) identities or (approximate) one-sided identities at least in the operator algebra case.
\begin{proposition}\label{pr:unique}
\begin{enumerate}
\item If $e$ is a quasi-identity of a ring, then so is $e^n$ for each $n\in\N$.
\item A contractive quasi-identity of a normed algebra is an idempotent, and hence its norm is either $0$ or $1$.
\item A contractive quasi-identity of an operator algebra $\A\subset\B(\Hi)$ is unique if it exists, and is Hermitian (hence an orthogonal projection).
\end{enumerate}
\end{proposition}
\begin{proof}$\underline{\text{(1):}}$ Let $e$ be a quasi-identity of a ring $\R$. For brevity of writing, we add an identity $1$ to $\R$ if it does not have one. Then $\left(\sum_{k=1}^n(1-e)e^{k-1}\right)r\left(\sum_{l=1}^ne^{l-1}(1-e)\right)=0,\forall r\in\R$ since $(1-e)\R(1-e)=\{0\}$. But each series is a ``telescoping series'', and the equation is simplified to $(1-e^n)r(1-e^n)=0,\forall r\in\R$, which means that $e^n$ is a quasi-identity of $\R$.

$\underline{\text{(2):}}$ Let $e$ be a contractive quasi-identity of a normed algebra. Then $e=2e^2-e^3$, that is, $e(e-e^2)=e-e^2$. Therefore inductively $e^n(e-e^2)=e-e^2,\forall n\in\N$. Thus $e-e^{n+1}=\sum_{k=1}^n(e^k-e^{k+1})=\sum_{k=1}^ne^{k-1}(e-e^2)=n(e-e^2)$. If $e^2\ne e$, then $1\ge\|e^{n+1}\|\ge n\|e-e^2\|-\|e\|,\forall n\in\N$. Therefore $e^2=e$.

$\underline{\text{(3):}}$ Let $e$ and $e'$ be contractive quasi-identities of an operator algebra $\A\subset\B(\Hi)$, hence they are idempotents by (2), and hence they are Hermitian by Lemma~\ref{lm: idempotent}~(1). Since $e$ is a quasi-identity, $e(e'-e'e)=e'-e'e$. Multiplying both sides by $e'$ on the left and right yields that $e'e(e'-e'ee')=e'-e'ee'$ since $e'$ is an idempotent. Therefore inductively $(e'e)^n(e'-e'ee')=e'-e'ee',\forall n\in\N$. Thus $e'-(e'e)^ne'=\sum_{k=1}^n(e'e)^{k-1}(e'-e'ee')=n(e'-e'ee')$, and hence $1\ge\|(e'e)^ne'\|\ge n\|e'-e'ee'\|-\|e'\|,\forall n\in\N$. Therefore
\begin{equation}\label{eq:e'=e'ee'}e'=e'ee',
\end{equation}and so $e'e=(e'e)^2$, that is, $e'e$ is an idempotent. By Lemma~\ref{lm: idempotent}~(1), $e'e=(e'e)^*=ee'$. Thus by Equation~(\ref{eq:e'=e'ee'}), $e'=(ee')e'=ee'$ since $e'$ is an idempotent. By symmetry, $e=e'e$, and hence $e'=ee'=e'e=e$.
\end{proof}In particular, if an operator algebra has a contractive one-sided or two-sided identity, then it is the only contractive quasi-identity.

As the following proposition shows that if a $C^*$-algebra has a quasi-identity (contractiveness is not assumed a priori), then it is necessarily an identity.
\begin{proposition}\label{pr:identity}If $\A$ is a $C^*$-algebra, then $\A$ possesses a quasi-identity if and only if $\A$ is unital. In this case, the identity is the only quasi-identity.
\end{proposition}
\begin{proof}Let $\A$ be a nonzero $C^*$-algebra, and let $\{e_{\alpha}\}$ be an approximate identity of $\A$. Suppose that $\A$ has a quasi-identity $e$. We may assume that $\A\subset\B(\Hi)$ nondegenerate, and denote the identity of $\B(\Hi)$ by 1. Then $(1-e)a(1-e)=0,\;\forall a\in\A$. In particular, for $\xi\in\Hi$, $(1-e)(e_{\alpha}-e)^*(1-e)\xi=0$. By taking the limit $\alpha\to\infty$, we have that $(1-e)(1-e)^*(1-e)\xi=0$. Since $\xi\in\Hi$ is arbitrary, $(1-e)(1-e)^*(1-e)=0$, so that $(1-e)^*(1-e)(1-e)^*(1-e)=0$, and thus $1=e\in\A$.
\end{proof}
\begin{corollary}\label{cor:ideal}If $J$ is a left (respectively, right) ideal in a $C^*$-algebra, then $J$ possesses a contractive quasi-identity if and only if $J$ has a contractive right (respectively, left) identity. In this case, the contractive right (respectively, left) identity is the only contractive quasi-identity.
\end{corollary}
\begin{proof}The second statement was already observed after Proposition~\ref{pr:unique}, and the ``if'' direction of the first statement is trivial. Assume that a left ideal $J$ in a $C^*$-algebra $\A\subset\B(\Hi)$ has a contractive quasi-identity $e$. Then $e$ is also a quasi-identity of the weak$^*$-closure of $\overline{J}^{\operatorname{w}^*}$ of $J$ in $\B(\Hi)$. Let $f$ be the identity of the von Neumann algebra $\overline{J^*J}^{\operatorname{w}^*}$ which is a subalgebra of $\overline{J}^{\operatorname{w}^*}$. Then $f$ is a contractive right identity of $\overline{J}^{\operatorname{w}^*}$. Obviously, $fef$ is a contractive quasi-identity of $\overline{J^*J}^{\operatorname{w}^*}$, and hence by Proposition~\ref{pr:identity}, $fef=f$. Thus $e=ef=(ef)^*=fe=fef=f$, where we used the fact that $e$ and $f$ are Hermitian.
\end{proof}
We are now in a position to present the characterization theorems.
\begin{theorem}\label{th:extinj}Let $X$ be a nonzero operator space, $z\in\Ball(\Q(X))$, and $(X,m_z)$ be the corresponding operator algebra.
\begin{enumerate}
\item$(X,m_z)$ has a quasi-identity of norm $1$ if $z\in\ext(\Ball(X^*))$.
\item$(X,m_z)$ has a left identity of norm $1$ if and only if $z\in X^*\cap\UL(\Q(X))$.
\item$(X,m_z)$ has a right identity of norm $1$ if and only if $z\in X^*\cap\UR(\Q(X))$.
\item$(X,m_z)$ has a two-sided identity of norm $1$ if and only if $z\in X^*\cap\U(\Q(X))$.
\end{enumerate}In each statement, $z^*$ is the quasi- (respectively, left, right, two-sided) identity of norm $1$.
\end{theorem}
\begin{proof}To see (1), assume that $z^*\in\ext(\Ball(X))$, then $\|z\|=1$. Let $\TER(X):=X\cap\Q(X)^*$ as in \cite{KP2}~Definition~4.6. Then $\TER(X)$ is a TRO, and $z^*\in\ext(\Ball(\TER(X)))$. Thus by Kadison's theorem (Lemma~\ref{lm: kadison}),
\begin{equation}\label{eq:kadison}(1_{11}-z^*z)\TER(X)(1_{22}-zz^*)=\{0\}.
\end{equation}Choosing $z^*\in\TER(X)$ yields that $(1_{11}-z^*z)z^*(1_{22}-zz^*)=0$, and so
$(1_{11}-z^*z)^2z^*=0$, and hence $z(1_{11}-z^*z)^2z^*=0$. Thus $z(1_{11}-z^*z)=0$, and $z^*z(1_{11}-z^*z)=0$.
Therefore, $z^*z$ is an idempotent.\footnote{We showed algebraically that $z^*z$ is an idempotent. Another way to see
this from Equation~(\ref{eq:kadison}) is to use spectral theory as in the proof of \cite{Kad}~Theorem~1 or
\cite{Pe}~Proposition~1.4.7, i.e., to consider the commutative $C^*$-algebra generated by $1_{11}$ and $z^*z$.} We
claim that $(1_{11}-z^*z)X(1_{22}-zz^*)=\{0\}$. Suppose the contrary, and pick $x_0\in X$ with $\|x_0\|\le1$ such that
$x_0=(1_{11}-z^*z)x_0(1_{22}-zz^*)\ne0$. Then $\|z^*\pm x_0\|^2=\|zz^*+x_0^*x_0\|=\max\{\|z\|^2,\|x_0\|^2\}=1$, and so
$z^*\pm x_0\in\Ball(X)$, and $z^*=\frac{1}{2}(z^*+x_0)+\frac{1}{2}(z^*-x_0)$. This contradicts the fact
that $z^*\in\ext(\Ball(X))$. Thus $(1_{11}-z^*z)X(1_{22}-zz^*)=\{0\}$ as claimed, i.e.,
$x=z^*zx+xzz^*-z^*zxzz^*\in X,\forall x\in X$, which tells us that $z^*\in X$ is a quasi-identity, and (1) has been shown. (3) was observed in \cite{KP2}~Proposition~2.10, and (2) is similar by symmetry, and (4) follows from (2) and (3).
\end{proof}
\begin{remark}\label{rm:extinj}
\begin{enumerate}
\item The converse direction in (1) does not hold. In fact, as we saw toward the beginning of this section, let $X:=\left[\begin{matrix}\C&\C&\C\\O&O&\C\\O&O&\C\end{matrix}\right]\subset\M_3(\C)$ with the usual matrix norm inherited from the operator norm of $\M_3(\C)$. Then $\Q(X)=\left[\begin{matrix}\C&\C&\C\\\C&\C&\C\\O&\C&\C\end{matrix}\right]$. Let $z=I_3\in\Q(X)$, where $I_3$ denotes the identity matrix. Then $(X,m_z)$ has a contractive quasi-identity $E_1+E_3$, but $z^*$ is not in $X$.
\item By Corollary~\ref{cor: subext}, $z^*$ in (2)--(4) is an extreme point of $\Ball(I(X))$, and hence an extreme point of $\Ball(\Q(X)^*)$, $\Ball(\T(X))$, $\Ball(X)$, $\Ball(\T(X)\cap\Q(X)^*)$ and $\Ball(\TER(X))$ too. However, in (1), $z^*$ is not an extreme point of $\Ball(\Q(X)^*)$ (hence not an extreme point of $\Ball(I(X))$) or $\Ball(\T(X))$ in general, though it is an extreme point of $\Ball(\TER(X))$ as stated in the proof above. In fact, let $X$ be as in the example above. Then $\T(X)=\M_3(\C)$. Let $z=E_1+E_3\in\Q(X)$. Then $z^*\in\ext(\Ball(X))$. But $z^*$ is not an extreme point of $\Ball(\Q(X)^*)$ or $\Ball(\T(X))$.
\item That $z\in\ext(\Ball(\Q(X)^*))$ does not imply that $(X,m_z)$ has an approximate quasi-identity. To see this, let $X:=\left[\begin{matrix}\C&\C&\C\\\C&\C&\C\\O&\C&\C\end{matrix}\right]$, which is a ``dual'' of the example above. Then $\Q(X)=\left[\begin{matrix}\C&\C&\C\\O&O&\C\\O&O&\C\end{matrix}\right]$. Let $z=E_1+E_3\in\Q(X)$. Then $z\in\ext(\Ball(\Q(X)))$, but $(X,m_z)$ does not have an approximate quasi-identity.
\end{enumerate}
\end{remark}The following corollary immediately follows from the theorem above.
\begin{corollary}\label{cor:extinj}Let $X$ be a nonzero operator space.
\begin{enumerate}
\item Some algebrization of $X$ has a quasi-identity of norm $1$ if and only if $\ext(\Ball(X))\cap\Q(X)^*\ne\varnothing$.
\item Some algebrization of $X$ has a left identity of norm $1$ if and only if $X\cap\UL(\Q(X))\ne\varnothing$
\item Some algebrization of $X$ has a right identity of norm $1$ if and only if $X\cap\UR(\Q(X))\ne\varnothing$.
\item Some algebrization of $X$ has a two-sided identity of norm $1$ if and only if $X\cap\U(\Q(X))\ne\varnothing$.
\end{enumerate}
\end{corollary}
\begin{proof}The only nontrivial part is the ``only if'' direction of (1). Let $e\in X$ with $\|e\|=1$ be a quasi-identity of $(X,m_z)$ for some $z\in\Ball(\Q(X))$. Then by Proposition~\ref{pr:unique}~(2), $eze=m_z(e,e)=e$ and so $ez$ and $ze$ are idempotents, and hence they are Hermitian by Lemma~\ref{lm: idempotent}~(1). Thus by the same argument as after Equation~(\ref{eq:ez=ezez,ze=zeze}) in the proof of Corollary~\ref{cor: inj}, we obtain that $e^*e=ze$. Similarly, $ee^*=ez$. Therefore $\forall x,y\in X, xe^*y=xe^*z^*e^*y=xe^*ezy=xzezy\in X$, and hence $e^*\in\Q(X)$. Since $e$ is a quasi-identity of $(X,m_z)$, $(1_{11}-ez)X(1_{22}-ze)=\{0\}$. But since $ez=ee^*$ and $ze=e^*e$, we have that $(1_{11}-ee^*)X(1_{22}-e^*e)=\{0\}$. Now we can show that $e\in\ext(\Ball(X))$ by the same way as to prove Kadison's theorem (Lemma~\ref{lm: kadison}), although $X$ is not a TRO in general. See the proof of \cite{Pe}~Proposition~1.4.7 for example.
\end{proof}Note that in the proof above, $z$ need not be $e^*$ in general. See the example in Remark~\ref{rm:extinj}~(1).

Extreme points best match quasi-identities when an operator space is injective as the following corollary shows. This fact convinces us that our attempt to characterize extreme points in terms of quasi-identities is a correct direction, and that defining multipliers (especially, quasi-multipliers) with the use of injective envelopes is the most plausible way. However, we remark that the alternative definitions $\Q''(X)$ and $\Q_{\pi}(X)$ defined in Section~\ref{section:alternative} are also useful in some occasions as we will see in Theorem~\ref{th:main}, and that $\Q(X)$ in Theorem~\ref{th:extinj}, Corollary~\ref{cor:extinj}, and Corollary~\ref{cor: inj} can be replaced by $\Q''(X)$ or $\Q_{\pi}(X)$ by the equivalence of the definitions.
\begin{corollary}\label{cor: inj}Let $X$ be a nonzero injective operator space, $z\in\Ball(\Q(X))$, and $(X,m_z)$ be the corresponding operator algebra. Then the following assertions hold.
\begin{enumerate}
\item$(X,m_z)$ has a quasi-identity of norm 1 if and only if $z\in\ext(\Ball(\Q(X)))$.
\item$(X,m_z)$ has a left identity of norm 1 if and only if $z\in\UL(\Q(X))$.
\item$(X,m_z)$ has a right identity of norm 1 if and only if $z\in\UR(\Q(X))$.
\item$(X,m_z)$ has a two-sided identity of norm 1 if and only if $z\in\U(\Q(X))$.
\end{enumerate}
\end{corollary}
\begin{proof}Note that if an operator space $X$ is injective, then $X=I(X)=\Q(X)^*$, and hence $z^*\in X$. Therefore all assertions follow from Theorem~\ref{th:extinj} except for the ``only if'' direction of (1). To show this direction, let $e$ be a quasi-identity of norm $1$. Then
\begin{equation}\label{eq:(1-ez)X(1-ze)=O}(1_{11}-ez)X(1_{22}-ze)=\{0\}.
\end{equation}By multiplying both sides by $z$ on the right, we have that $(1_{11}-ez)Xz(1_{11}-ez)=\{0\}$. By choosing $(1_{11}-ez)^*z^*\in X$ and multiplying both sides by $(1_{11}-ez)^*$ on the right, we have that $(1_{11}-ez)(1_{11}-ez)^*z^*z(1_{11}-ez)(1_{11}-ez)^*=0$, which implies that $z(1_{11}-ez)(1_{11}-ez)^*=0$. Hence $z(1_{11}-ez)(1_{11}-ez)^*z^*=0$, and accordingly, $z(1_{11}-ez)=0$.
\begin{equation}\label{eq:z=zez}z=zez.
\end{equation}Thus
\begin{equation}\label{eq:ez=ezez,ze=zeze}ez=ezez\quad\text{and}\quad ze=zeze,
\end{equation}which means that $ez$ and $ze$ are idempotents.\footnote{Equations~(\ref{eq:ez=ezez,ze=zeze}) is actually an immediate consequence of Proposition~\ref{pr:unique}~(2). But we deduced Equation~(\ref{eq:z=zez}) since we use it toward the end of the proof.} Therefore by Lemma~\ref{lm: idempotent}~(1), $ez$ and $ze$ are orthogonal projections, and hence
\begin{equation}\label{eq:ez=z*e*,ze=e*z*}ez=(ez)^*=z^*e^*\text{\quad and\quad}ze=(ze)^*=e^*z^*.
\end{equation}Hence Equation~(\ref{eq:(1-ez)X(1-ze)=O}) is rewritten to $(1_{11}-z^*e^*)X(1_{22}-e^*z^*)=\{0\}$, and by repeating\footnote{Alternatively, one can obtain Equation~(\ref{eq:e*=e*z*e*}) from Proposition~\ref{pr:unique}~(2). In fact, $e=m_z(e,e)=eze$.} the argument above Equation~(\ref{eq:z=zez}), we obtain that
\begin{equation}\label{eq:e*=e*z*e*}e^*=e^*z^*e^*.
\end{equation}Thus $e^*e\ge e^*z^*ze\ge e^*z^*e^*eze=e^*e$, and together with Equations~(\ref{eq:ez=z*e*,ze=e*z*})~and~(\ref{eq:ez=ezez,ze=zeze}), we have that $e^*e=e^*z^*ze=zeze=ze$. Hence together with Equations~(\ref{eq:e*=e*z*e*}),~(\ref{eq:ez=z*e*,ze=e*z*}),~and~(\ref{eq:z=zez}), we have that $e^*=e^*z^*e^*=e^*ez=zez=z$. Thus Equation~(\ref{eq:(1-ez)X(1-ze)=O}) becomes $(1_{11}-z^*z)X(1_{22}-zz^*)=\{0\}$, which tells that $z^*\in\ext(\Ball(X))$ by Kadison's theorem (Lemma~\ref{lm: kadison}).
\end{proof}The alternative definitions of multipliers which we defined in Section~\ref{section:alternative} work in the ``approximate'' version of characterization.
\begin{theorem}\label{th:main}Let $X$ be an operator space, $z$ be in $\Ball(\Q''(X))$ or $\Ball(\Q_{\pi}(X))$, and $(X,m_z)$ be the corresponding operator algebra. Then the following implications hold.

\begin{tabular}{rcrl}\em{(1)}&&\em{(i)}&$z^*\in\ext(\Ball\left(\overline{X}^{
\operatorname{w}^*}\right))$, where the weak$^*$-closure is taken in $C^*(\partial X)''$, and\\&&&$z^*z\in\LM''(X)$ and $zz^*\in\RM''(X)$;\\&$\Rightarrow$&\em{(ii)}&$(X,m_z)$ has a contractive approximate quasi-identity;\\&$\Leftarrow$&\em{(iii)}&$z^*\in\ext(\Ball\left(\overline{X}
^{\operatorname{w}^*}\right))$, where the weak$^*$-closure is taken in $\B(\Hi_1\oplus\Hi_2)$, and\\&&&$z^*z\in\LM_{\pi}(X)$ and $zz^*\in\RM_{\pi}(X)$.\\
\em{(2)}&&\em{(i)}&$z\in\overline{X^*}^{\operatorname{w}^*}\cap\UL(\Q''(X))$, where the weak$^*$-closure is taken in $C^*(\partial X)''$;\\
&$\Leftrightarrow$&\em{(ii)}&$(X,m_z)$ has a contractive approximate left identity;\\
&$\Leftrightarrow$&\em{(iii)}&$z\in\overline{X^*}^{\operatorname{w}^*}\cap\UL(\Q_{\pi}(X))$, where the weak$^*$-closure is taken in $\B(\Hi_1\oplus\Hi_2)$.\\
\em{(3)}&&\em{(i)}&$z\in\overline{X^*}^{\operatorname{w}^*}\cap\UR(\Q''(X))$, where the weak$^*$-closure is taken in $C^*(\partial X)''$;\\
&$\Leftrightarrow$&\em{(ii)}&$(X,m_z)$ has a contractive approximate right identity;\\
&$\Leftrightarrow$&\em{(iii)}&$z\in\overline{X^*}^{\operatorname{w}^*}\cap\UR(\Q_{\pi}(X))$, where the weak$^*$-closure is taken in $\B(\Hi_1\oplus\Hi_2)$.\\
\em{(4)}&&\em{(i)}&$z\in\overline{X^*}^{\operatorname{w}^*}\cap\U(\Q''(X))$, where the weak$^*$-closure is taken in $C^*(\partial X)''$;\\
&$\Leftrightarrow$&\em{(ii)}&$(X,m_z)$ has a contractive approximate two-sided identity;\\
&$\Leftrightarrow$&\em{(iii)}&$z\in\overline{X^*}^{\operatorname{w}^*}\cap\U(\Q_{\pi}(X))$, where the weak$^*$-closure is taken in $\B(\Hi_1\oplus\Hi_2)$.
\end{tabular}
\end{theorem}
\begin{proof}We will show ``(i)$\Rightarrow$(ii)'' of (1). ``(iii)$\Rightarrow$(ii)'' is exactly the same. In the proof, all the weak$^*$-closures are taken in $C^*(\partial X)''$. Assume that $z^*\in\ext(\Ball\left(\overline{X}^{\operatorname{w}^*}\right))$. Then $\|z\|=1$, and $z^*$ is an extreme point of the unit ball of the weak$^*$-closed TRO $\WTER(X):=\overline{X}^{\operatorname{w}^*}\cap\overline{\Q''(X)^*}^{\operatorname{w}^*}$ as well. Thus by Kadison's theorem (Lemma~\ref{lm: kadison}), $(1_{\E}-z^*z)\WTER(X)
(1_{\F}-zz^*)=\{0\}$. Then by the same argument as in the proof of Theorem~\ref{th:extinj}~(1), we obtain that
\begin{equation}\label{eq:quasi}x=z^*zx+xzz^*-z^*zxzz^*,\quad\forall x\in\overline{X}^{\operatorname{w}^*}.
\end{equation}Pick a net $\{e_{\alpha}\}\subset X$ of contractions such that $\operatorname{w}^*$-$\lim_\alpha e_\alpha=z^*$. By the separate weak$^*$-continuity of the product in $C^*(\partial X)''$, for each $x\in X$, $\operatorname{w}^*$-$\lim_\alpha e_\alpha zx=z^*zx\in X$ since $z^*z\in\LM''(X)$. Thus
\begin{equation}\label{eq:left}\operatorname{w}\text{-}\lim_{\alpha\to\infty}e_\alpha zx=z^*zx\in X,\quad\forall x\in X.
\end{equation}Similarly,
\begin{equation}\label{eq:right}\operatorname{w}\text{-}\lim_{\alpha\to\infty}xze_\alpha=xzz^*\in X,\quad\forall x\in X.
\end{equation}Now we adopt a technique employed in the proof of Theorem~2.2 in \cite{ER1}. Let $\F$ be the collection of the finite subsets of $X$, and let $\Lambda:=\F\times\N$. Then $\Lambda$ is a directed set under the ordering ``$(F_1,n_1)\le(F_2,n_2)$ if and only if $F_1\subset F_2$ and $n_1\le n_2$''. Given $F=\{x_1,\dots,x_m\}\in\F$, let$$V_F:=\left\{\{(e-z^*)zx_1,\dots,(e-z^*)zx_m,x_1z(e-z^*),\dots,x_mz(e-z^*)\}\;;\;e\in \Ball(X)\right\}\subset X^{2m},$$where $X^{2m}$ is given the supremum norm. It follws from Equations~(\ref{eq:left})~and~(\ref{eq:right}) that $\vec{0}:=(0_1,\dots,0_{2m})$ lies in the weak-closure of $V_F$ in $X^{2m}$, and hence it lies in the norm-closure of $V_F$ in $X^{2m}$ since $V_F$ is convex. Therefore for given $n\in\N$, $V_F\cap\{\vec{x}\in X^{2m};\|\vec{x}\|<1/n\}\ne\varnothing$. The argument above tells us that for given $(F,n)\in\Lambda$, we may choose $e_\lambda\in\Ball(X)$ with $\|(e_\lambda-z^*)zx_k\|<1/n$ and $\|x_kz(e_\lambda-z^*)\|<1/n$ for $k=1,\dots,m$. Hence we have obtained a contractive net $\{e_\lambda\}$ such that
\begin{equation}\label{eq: 3}\lim_{\lambda\to\infty}e_\lambda zx=z^*zx\in X,\quad\forall x\in X,
\end{equation}
\begin{equation}\label{eq: 4}\lim_{\lambda\to\infty}xze_\lambda=xzz^*\in X,\quad\forall x\in X,
\end{equation}By putting Equations~(\ref{eq:quasi}),~(\ref{eq: 3}),~and~(\ref{eq: 4}) altogether, we obtain that $x=\lim_\lambda\lim_{\lambda'}(e_\lambda zx+xze_{\lambda'}-e_\lambda zxze_{\lambda'}),\forall x\in X$. Since $\{e_{\lambda}\}$ is bounded, by the routine argument using the triangular inequality, we have that $x=\lim_\lambda(e_\lambda zx+xze_{\lambda}-e_\lambda zxze_{\lambda}),\forall x\in X$.

``(i)$\Rightarrow$(ii)$\Leftarrow$(iii)'' in (2)--(4) are similar but easier. To see ``(i)$\Leftarrow$(ii)$\Rightarrow$(iii)'' of (2)--(4), we show (ii)$\Rightarrow$(i) of (2). The others are similar. Let $\{e_\alpha\}\subset X$ be a contractive approximate left identity of $(X,m_z)$, and let $e$ be its weak$^*$ accumulation point in $\overline{X}^{\operatorname{w}^*}$. Then $(1_{\E}-ez)X=\{0\}$, so $ez=1_{\E}$. Thus by Lemma~\ref{lm: idempotent}~(2), $e=z^*$.
\end{proof}
\begin{remark}\label{rm:main}
\begin{enumerate}
\item The converse directions in (1) does not hold. See the example in Remark~\ref{rm:extinj}~(1).
\item By Corollary~\ref{cor:subext}~(1), $z^*$ in (i) of (2)--(4) is an extreme point of $\Ball\left(\T(X)''\right)$, and hence an extreme point of $\Ball\left(\overline{\Q''(X)^*}^{\operatorname{w}^*}\right)$, $\Ball(\Q''(X)^*)$, $\Ball\left(\overline{X}^{\operatorname{w}^*}\right)$, $\Ball(\WTER(X))$, and $\Ball\left(\overline{X}^{\operatorname{w}^*}\cap\Q''(X)^*\right)$ too. However, in (i) of (1), $z^*$ is not an extreme point of $\Ball(\Q''(X)^*)$ (hence not an extreme point of $\Ball\left(\overline{\Q''(X)^*}^{\operatorname{w}^*}\right)$ or $\Ball\left(\T(X)''\right)$) in general, though it is an extreme point of $\Ball(\WTER(X))$ (and hence an extreme point of $\Ball\left(\overline{X}^{\operatorname{w}^*}\cap\Q''(X)^*\right)$) as stated in the proof above. See the example in Remark~\ref{rm:extinj}~(1). Here all the weak$^*$-closures are taken in $C^*(\partial X)''$.
\item By Corollary~\ref{cor:subext}~(2), $z^*$ in (iii) of (2)--(4) is an extreme point of $\Ball(\K,\Hi)$, and hence an extreme point of $\Ball\left(\overline{\Q_{\pi}(X)^*}^{\operatorname{w}^*}\right)$, $\Ball(\Q_{\pi}(X))$, $\Ball\left(\overline{\T(X)^*}^{\operatorname{w}^*}\right)$, $\Ball\left(\overline{X}^{\operatorname{w}^*}\right)$, $\Ball\left(\overline{\T(X)^*}^{\operatorname{w}^*}\cap\overline{\Q_{\pi}(X)^*}^{\operatorname{w}^*}\right)$, $\Ball\left(\overline{\T(X)^*}^{\operatorname{w}^*}\cap\Q_{\pi}(X)^*\right)$, $\Ball(\WTER(X))$ \footnote{Here $\WTER(X):=\overline{X}^{\operatorname{w}^*}\cap\overline{\Q_{\pi}(X)^*}^{\operatorname{w}^*}$ with the weak$^*$-closures taken in $\B(\Hi_1\oplus\Hi_2)$.}, and $\Ball\left(\overline{X^*}^{\operatorname{w}^*}\cap\Q_{\pi}(X)^*\right)$ too. However, in (iii) of (1), $z^*$ is not an extreme point of $\Ball(\Q_{\pi}(X)^*)$ (hence not an extreme point of $\Ball\left(\overline{\Q_{\pi}(X)^*}^{\operatorname{w}^*}\right)$ or $\Ball\left(\overline{\T(X)^*}^{\operatorname{w}^*}\right)$) in general, though it is an extreme point of $\Ball(\WTER(X))$ and hence an extreme point of $\Ball\left(\overline{X^*}^{\operatorname{w}^*}\cap\Q_{\pi}(X)^*\right)$. See the example in Remark~\ref{rm:extinj}~(1). Here all the weak$^*$-closures are taken in $\B(\Hi_1\oplus\Hi_2)$.
\end{enumerate}
\end{remark}We imposed the condition ``$z^*z\in\LM''(X)$ and $zz^*\in\RM''(X)$'' or ``$z^*z\in\LM_{\pi}(X)$ and $zz^*\in\RM_{\pi}(X)$'' in (1) of the Theorem~\ref{th:main} to avoid the difficulty which would come from the fact that an operator algebra product is weak$^*$-continuous with respect to each factor ``separately''. This condition is so strong that even an operator algebra with a contractive approximate one-sided identity need not satisfy it. In fact, if a contractive approximate left identity $\{e_\alpha\}$ of an operator algebra $\A$ satisfies this condition, then Equation~(\ref{eq: 4}) in the proof above suggests that $\{e_\alpha\}$ can be chosen so that $\lim_\alpha ae_\alpha$ exists in $\A$ for all $a\in\A$. However, this is not possible in general as the following example shows. We thank David~P.~Blecher for the basic idea of the example.
\begin{example}\label{ex:B}Let us canonically identify $\B\left(\bigoplus_{n=1}^{\infty}l^2(\N)\right)$ with a subset of the set $\M(\B(l^2(\N)))$ of $\aleph_0\times\aleph_0$ matrices with entries in $\B(l^2(\N))$, and $\B(l^2(\N))$ with a subset of the set $\M(\C)$ of $\aleph_0\times\aleph_0$ matrices with entries in $\C$. Let$$\A:=\left\{[a_{i,j}]\in\B\left(\bigoplus_{n=1}^{\infty}l^2(\N)\right)\subset\M\left(\B\left(l^2(\N)\right)\right);
\exists i_0\in\N\text{ such that }a_{i,j}=0,\forall i\ge i_0,\forall j\in\N\right\}^-.$$Then $\A$ is a right ideal of the von Neumann algebra $\B\left(\bigoplus_{n=1}^{\infty}l^2(\N)\right)$, so it has a contractive approximate left identity.\footnote{For example, the sequence $\{e_n\}$ with $e_n\in\A$ whose first $n$ diagonal elements are the identity operators on $l^2(\N)$ and all the other elements are zero operators, is a contractive approximate left identity of $\A$.} For each $j\in\N$, let us denote by $E_j$ the element of $\B(l^2(\N))\subset\M(\C)$ whose $(j,j)$-entry is $1$ and all other entries are $0$'s. Define $a=[a_{i,j}]\in\A$ as follows: $a_{1,j}:=E_j,\forall j\in\N$; $a_{i,j}=0,\forall i\ge2,
\forall j\in\N$. Also for each $i\in\N$, define $p_i\in\M\left(\B\left(l^2(\N)\right)\right)$ as follows: The $(i,i)$-entry of $p_i$ is the identity operator on $\B\left(l^2(\N)\right)$, and all other entries of $p_i$ are zero operators. Then note that both $ap_i$ and $(ap_i)^*$ are in $\A$. Let $\{e_{\alpha}\}$ be ``any'' contractive approximate left identity of $\A$. Since $\A$ is separable, one can take a subnet of $\{e_{\alpha}\}$ to be a sequence $\{e_{\alpha_n}\}$. By \cite{B}~Lemma~2.2~(1), $\lim_ne_{\alpha_n}^*b=b,\forall b\in\A$. In particular, $\lim_ne_{\alpha_n}^*(ap_i)^*=(ap_i)^*$, and hence $\lim_nap_ie_{\alpha_n}=ap_i$. So one may choose the sequence above inductively in the following way: For each $n\in\N$, pick $m_n\in\N$ such that the entries of the $m$-th row of the $(1,j)$-entry of $ae_{\alpha_n}$ are all $0$'s for all $j\in\N$ and all $m\ge m_n$; $\|ap_{m_n}e_{\alpha_{n+1}}\|>1/2$. Now it is easy to see that $\|ae_{\alpha_{n+1}}-ae_{\alpha_n}\|>1/2,\forall n\in\N$, and so $\{ae_{\alpha_n}\}$ is not a Cauchy sequence. Therefore, $\lim_nae_{\alpha_n}$, and hence $\lim_{\alpha}ae_{\alpha}$ does not exist.
\end{example}The example above also yields the following proposition.
\begin{proposition}\label{pr:nondeg}There exists an operator space $X$ for which there is a nondegenerate representation of the $C^*$-algebra $I(\Sy_X)$ on a Hilbert space $\Hi$ such that the $C^*$-subalgebra $C^*(\partial X)$ is degenerate on $\Hi$.
\end{proposition}
\begin{proof}Let $X$ be the underlying operator space of the operator algebra $\A$ defined in Example~\ref{ex:B}, $z\in\Ball(\Q(X))$ be the quasi-multiplier associated with $\A$, and $\{e_{\alpha}\}$ be a contractive approximate left identity of $\A=(X,m_z)$. Represent the $W^*$-algebra $I(\Sy_X)''$ weak$^*$-continuously on a Hilbert space $\Hi$ nondegenerately. Let $e$ be a weak$^*$ accumulation point of $\{e_{\alpha}\}$ in $I(\Sy_X)''$. Then $x=\lim_{\alpha}e_{\alpha}zx=ezx,\forall x\in X$, and hence $(1_{11}-ez)\overline{X}^{\operatorname{w}^*}=\{0\}$. Suppose that $C^*(\partial X)$ is nondegenerate on $\Hi$. Then $1_{11}-ez=0$, and so $e=z^*$ by Lemma~\ref{lm: idempotent}~(2). Thus $xzz^*\in\overline{X}^{\operatorname{w}^*}\cap I(\Sy_X)=X$. Now one can choose a contractive approximate left identity $\{e_{\lambda}\}$ of $\A=(X,m_z)$ so that $\lim_{\lambda}m_z(x,e_{\lambda})=\lim_{\lambda}xze_{\lambda}\in X,\forall x\in X$ as in the proof of (i)$\Rightarrow$(ii) of (1) of Theorem~\ref{th:main}. This contradicts the fact observed in Example~\ref{ex:B}.
\end{proof}The question ``If an operator space $X$ is injective, then always $\ext(\Ball(X))\ne\varnothing$?'' naturally arises. If the answer is yes, then by Corollary~\ref{cor: inj}~(1), any injective operator space can be made into an operator algebra with a quasi-identity with norm 1. More generally, one may ask ``For any operator space $X$, $\ext(\Ball(\Q(X)))\ne\varnothing$?'' The answer to the second question being yes yields the answer to the first question being yes since $X=\Q(X)^*$ for an injective operator space $X$. We leave them as well as the following as open questions. Question: Can a TRO always be made into an operator algebra with a contractive approximate quasi-identity?

We close this section by recalling two examples from \cite{KP2}. The quasi-multiplier space of the operator space $X$
in Example~2.13 of that paper is $\{0\}$. So the ``zero product'' is the only possible operator algebra product that
$X$ can be equipped with. And hence, there is no algebrization for $X$ to have a quasi-identity. A more interesting
example is Example~2.11 of the same paper. Elementary but tedious calculations show that all points on the sphere
(i.e., the set of points with norm 1) of $\mathcal{X}$ and $\Q(\mathcal{X})$ are extreme points, and $\mathcal{X}$ can have a quasi-identity for a certain algebrization, however, there is no algebrization for $\mathcal{X}$ to have a
``contractive'' quasi-identity.
\section{$C^*$-algebras and their one-sided ideals}\label{section: C*}
In this section, we give an operator space characterization of $C^*$-algebras and their one-sided ideals in terms of quasi-multipliers. Although we prefer to use $\Q(X)$, the reader should keep in mind that these characterizations can be formulated using alternative definitions $\Q''(X)$ or $\Q_{\pi}(X)$.

First we characterize one-sided ideals in $C^*$-algebras. Another characterization of such ideals was given on page 2108 of \cite{BK}: Left ideals in $C^*$-algebras are exactly the operator algebras $A$ with a r.c.a.i. that are also abstract triple systems.
\begin{theorem}\label{th:ideal}Let $X$ be a nonzero operator space, and $z\in\Ball(\Q(X))$, and $(X,m_z)$ be the corresponding operator algebra. Then there is a completely isometric homomorphism from $(X,m_z)$ onto a left (respectively, right) ideal in some $C^*$-algebra if and only if $z\in\UR(\Q(X))$, $X^*X\subset zX$, and $Xz\subset XX^*$ (respectively, $z\in\UL(\Q(X))$, $XX^*\subset Xz$, and $zX\subset X^*X$).
\end{theorem}
\begin{proof}We prove the left ideal case. The right ideal case is similar by symmetry.

\underline{$\Rightarrow$:} Assume that $(X,m_z)$ is a left ideal in a $C^*$-algebra. Then it has a contractive approximate right identity $\{e_{\alpha}\}$. As in the proof of \cite{BK}~Theorem~2.3, there is a $v\in I(X)$ such that $xv^*=\lim_{\alpha}xe_{\alpha}^*\in XX^*,\forall x\in X$ and $v^*v=1_{22}$, where the products are taken in the injective $C^*$-algebra $I(\Sy_X)$. Thus $v^*X=v^*XX^*X\supset\lim_{\alpha}v^*e_{\alpha}X^*X=v^*vX^*X=X^*X$. By the first sentence in the proof of Lemma~\ref{lm:BK} of the present paper, $z=v^*$. Hence $Xz\subset XX^*$, $z\in\UR(\Q(X))$, and $X^*X\subset zX$.

\underline{$\Leftarrow$:} That $X^*X\subset zX$ implies that $XX^*X=XzX\subset X$ since $z\in\Q(X)$. So $X$ is a TRO, and hence $XX^*$ is a $C^*$-algebra. Define $\psi:X\to XX^*$ by $\psi(x):=xz,\forall x\in X$. Then $\psi$ is a completely contractive mapping from $X$ into the $C^*$-algebra $XX^*$. In fact $\psi$ is a complete isometry since the right multiplication by the contractive element $z^*$ gives the inverse mapping of $\psi$, i.e., $\psi(x)z^*=xzz^*=x,\forall x\in X$ since $z\in\UR(\Q(X))$. That $\psi(m_z(x_1,x_2))=x_1zx_2z=\psi(x_1)\psi(x_2),\forall x_1,x_2\in X$ shows that $\psi$ is a homomorphism. Since $XX^*Xz\subset Xz$, $\psi(X)=Xz$ is a left ideal in the $C^*$-algebra $XX^*$.
\end{proof}Now we give an operator space characterization of $C^*$-algebras in terms of quasi-multipliers. It makes a beautiful contrast with the one-sided ideal case above. We remark that the ``$\Leftarrow$'' directions of the theorem below was essentially first observed by Vern~I.~Paulsen assuming that $X$ is a unital $C^*$-algebra and using the classical definition of quasi-multipliers (\cite{Pe}~Section~3.12), in which case $\Q(X)=X$. We thank him for letting us know his observation. In the following theorem and its proof we revive the symbols $\odot$ and $\bullet$ defined in Section~\ref{section: pre} to avoid confusion.
\begin{theorem}\label{th: C*} Let $X$ be a nonzero operator space, and $z\in\Ball(\Q(X))$, and $(X,m_z)$ be the corresponding operator algebra. Then $(X,m_z)$ is a $C^*$-algebra with a certain involution $\sharp$ if and only if $z\in\U(\Q(X))$ and $X\bullet z=z^*\bullet X^*$ (or, equivalently\footnote{These equivalences are obvious since $z^*\bullet z=1_{11}$ and $z\bullet z^*=1_{22}$.}, $z\bullet X=X^*\bullet z^*$, or $z^*\bullet X^*\bullet z^*=X$). The involution $\sharp$ is uniquely given by $x^\sharp=z^*\bullet x^*\bullet z^*,\;\forall x\in X$, which implies that for a given operator algebra there exists at most one involution that makes the operator algebra a $C^*$-algebra. Moreover, all such $C^*$-algebras are $*$-isomorphic, which recovers the no doubt well-known fact that for a given operator space there exists at most one $C^*$-algebra structure up to $*$-isomorphism.
\end{theorem}
\begin{proof}$\underline{\Leftarrow:}$ Let $z\in\U(\Q(X))$ and $z\bullet X=X^*\bullet z^*$. Define an involution $\sharp$ by $x^{\sharp}:=z^*\bullet x^*\bullet z^*,\;\forall x\in X$, where $*$ is the involution on the injective envelope $C^*$-algebra $I(\Sy_X)$. Since $z\bullet X=X^*\bullet z^*$, $z^*\bullet x^*\bullet z^*$ is certainly in $X$. And also $(x^\sharp)^\sharp=z^*\bullet z\bullet x\bullet z\bullet z^*=1_{11}\bullet x \bullet 1_{22}=x$. Hence $\sharp$ is a well-defined involution. $\|m_z(x^\sharp,x)\|=\|z^*\bullet x^*\bullet z^*\bullet z\bullet x\|=\|z^*\bullet x^*\bullet 1_{11}\bullet x\|=\|z^*\bullet x^*\bullet x\|\ge\|z\bullet z^*\bullet x^*\bullet x\|=\|1_{22}\bullet x^*\bullet x\|=\|x^*\bullet x\|=\|x\|^2$ shows that $(X,m_z,\sharp)$ is a $C^*$-algebra.

$\underline{\Rightarrow:}$ Assume that $(X,m_z,\sharp)$ is a $C^*$-algebra. By Theorem~\ref{th:extinj}~(4), $z\in\U(\Q(X))$. To check that $X\bullet z=z^*\bullet X^*$, we may assume that $(X,m_z,\sharp)\subset\B(\K)$ as a $C^*$-subalgebra for some Hilbert space $\K$. Let $\Sy_X':=\left[\begin{matrix}\C1_{\K}&X\\X^\sharp&\C1_{\K}\end{matrix}\right]\subset\M_2(\B(\K))\;(\text{Actually $X^\sharp=X$.})$ be Paulsen's operator system, and $C^*(X)=\M_2(X)$ be the $C^*$-algebra generated by $\left[\begin{matrix}O&X\\O&O\end{matrix}\right]$ in $\M_2(\B(\K))$. By using Hamana's theorem (\cite{H3}~Corollary~4.2) it is easily seen that there is a $*$-homomorphism $\Psi=\left[\begin{matrix}\Psi_{11}&\Psi_{12}\\\Psi_{21}&\Psi_{22}\end{matrix}\right]$, which is factored by a well-known trick, from $C^*(X)$ onto $C^*(\partial X)$ such that $\Psi_{12}(x)=x$ (and hence $\Psi_{21}(x^\sharp)=(\Psi_{12}(x))^*$), $\forall x\in X$, where $C^*(\partial X)$ is as in Section~\ref{section: pre}.
Let $\{e_\alpha\}$ be a contractive approximate identity of the $C^*$-algebra $(X,m_z,\sharp)$. Then $\left[\begin{matrix}0&m_z(x,y)\\0&0\end{matrix}\right]=\lim_\alpha\Psi\left(\left[\begin{matrix}0&x\\0&0\end{matrix}
\right]\left[\begin{matrix}0&0\\e_\alpha^\sharp&0\end{matrix}\right]\left[\begin{matrix}0&y\\0&0\end{matrix}\right]
\right)=\lim_\alpha\left[\begin{matrix}0&x\\0&0\end{matrix}\right]\odot\left[\begin{matrix}0&0\\e_\alpha^*&0\end{matrix}
\right]\odot\left[\begin{matrix}0&y\\0&0\end{matrix}\right]=\lim_\alpha\left[\begin{matrix}0&x\bullet e_\alpha^*\bullet y\\0&0\end{matrix}\right],\;\forall x,y\in X$, so that $\lim_\alpha x\bullet e_\alpha^*\bullet y=x\bullet z\bullet y,\;\forall x,y\in X$. Now $\lim_\alpha\left[\begin{matrix}x\bullet e_\alpha^*&0\\0&0\end{matrix}\right]=\lim_\alpha\Psi\left(\left[\begin{matrix}0&x\\0&0\end{matrix}\right]\right)\odot
\Psi\left(\left[\begin{matrix}0&0\\e_\alpha^\sharp&0\end{matrix}\right]\right)=\lim_\alpha\Psi\left(\left[\begin{matrix}
xe_\alpha^\sharp&0\\0&0\end{matrix}\right]\right)=\left[\begin{matrix}\Psi_{11}(x)&0\\0&0\end{matrix}\right],\;\forall x\in X$. Thus $\Psi_{11}(x)\bullet y=x\bullet z\bullet y,\;\forall x,y\in X$, and hence by Lemma~\ref{lm: BP}~(1), $\Psi_{11}(x)=x\bullet z,\;\forall x\in X$, so that $\Psi_{11}(X)=X\bullet z$. On the other hand, $\lim_\alpha\left[\begin{matrix}e_\alpha\bullet x^*&0\\0&0\end{matrix}\right]=\lim_\alpha\Psi\left(\left[\begin{matrix}0&e_\alpha\\0&0\end{matrix}\right]\right)\odot
\Psi\left(\left[\begin{matrix}0&0\\x^\sharp&0\end{matrix}\right]\right)=\lim_\alpha\Psi\left(\left[\begin{matrix}
e_\alpha x^\sharp&0\\0&0\end{matrix}\right]\right)=\left[\begin{matrix}\Psi_{11}(x^\sharp)&0\\0&0\end{matrix}\right],\;\forall x\in X$. Thus $y^*\bullet\Psi_{11}(x^\sharp)=y^*\bullet z^*\bullet x^*,\;\forall x,y\in X$, since $\lim_\alpha y^*\bullet e_\alpha\bullet x^*=y^*\bullet z^*\bullet x^*,\;\forall x,y\in X$. Hence by Lemma~\ref{lm: BP}~(1) again,
$\Psi_{11}(x^\sharp)=z^*\bullet x^*,\;\forall x\in X$, so that $\Psi_{11}(X)=z^*\bullet X^*$. Therefore, $X\bullet z=\Psi_{11}(X)=z^*\bullet X^*$. It also follows that $z^*\bullet x^*\bullet z^*=\Psi_{11}(x^\sharp)\bullet z^*=x^\sharp\bullet z\bullet z^*=x^\sharp,\;\forall x\in X$.

Finally, we show that all $C^*$-algebras which have the same underlying operator space $X$ are $*$-isomorphic. Since this fact is no doubt well known, and a simpler proof (or observation) is possible, it might be redundant to present the proof. However, it would be instructive to show how two quasi-multipliers work out, so we include the proof. Let $z'\in\U(\Q(X))$, and assume that $(X,m_{z'},\natural)$ is also a $C^*$-algebra. Then the involution $\natural$ is given by $x^\natural=z'^*\bullet x^*\bullet z'^*$. Define a linear mapping $\pi:(X,m_z,\sharp)\to(X,m_{z'},\natural)$ by $x\mapsto x\bullet z\bullet z'^*$. We must check that the image is certainly in $X$. Note that $\Psi_{11}(x)=x\bullet z$ which is one-to-one, and $(\Psi_{11})^{-1}(a)=a\bullet z^*,\;\forall a\in\Psi_{11}(X)$. Considering $(X,m_{z'},\natural)\subset\B(\K')$ as a $C^*$-subalgebra for some Hilbert space $\K'$, we can define $\Psi_{11}'$ as we defined $\Psi_{11}$. Then $(\Psi_{11}')^{-1}(a)=a\bullet z'^*,\;\forall a\in\Psi_{11}'(X)$. By noting that $\Psi_{11}(X)=\Psi_{11}'(X)=\E(X)$, where $\E(X)$ is as in Section~\ref{section: pre}, we have that for $x\in X,\;x\bullet z\bullet z'^*=(\Psi_{11}')^{-1}(\Psi_{11}(x))\in X$, so that Im$\pi\subset X$. Similarly, we have that $x\bullet z'\bullet z^*=(\Psi_{11})^{-1}(\Psi_{11}'(x))\in X,\;\forall x\in X$. Thus $x=x\bullet z'\bullet z^*\bullet z\bullet z'^*=\pi(x\bullet z'\bullet z^*),\;\forall x\in X$, which shows that $\pi$ is onto. $\pi$ being one-to-one follows from $x=\pi(x)\bullet z'\bullet z^*,\;\forall x\in X$. Furthermore, $\pi(m_z(x,y))=x\bullet z\bullet y\bullet z\bullet z'^*=x\bullet z\bullet z'^*\bullet z'\bullet y\bullet z\bullet z'^*=m_{z'}(\pi(x),\pi(y))$ and $\pi(x^\sharp)=z^*\bullet x^*\bullet z^*\bullet z\bullet z'^*=z^*\bullet x^*\bullet z'^*=z'^*\bullet z'\bullet z^*\bullet x^*\bullet z'^*=z'^*\bullet(x\bullet z\bullet z'^*)^*\bullet z'^*=\pi(x)^\natural$ show that $\pi:(X,m_z,\sharp)\to(X,m_{z'},\natural)$ is a $*$-homomorphism.
\end{proof}One may expect that the quasi-multiplier space of a $C^*$-algebra always can be a $C^*$-algebra for some algebrization, or the quasi-multiplier space of a TRO is a TRO. However, neither of them is true in general. The following example shows that the quasi-multiplier space of a $C^*$-algebra may not even be a TRO, hence may not be completely isometric to a one-sided ideal in any $C^*$-algebra.
\begin{example}\label{ex: counter}Let $\Hi$ be an infinite dimensional Hilbert space and let $\Kp(\Hi)^1$ denote the unitization of $\Kp(\Hi)$ by the identity $1$ of $\B(\Hi)$, where $\Kp(\Hi)$ is the set of the compact operators on $\Hi$. Define $X:=\left[\begin{matrix}\Kp(\Hi)&\Kp(\Hi)\\\Kp(\Hi)&\Kp(\Hi)^1\end{matrix}\right]$. Give $X$ the canonical operator space structure as a subspace of $\B(\Hi\oplus\Hi)$, then $X$ is a $C^*$-algebra with the product on $\B(\Hi\oplus\Hi)$. It is easy to see that the product $\bullet$ defined in Section~\ref{section: pre} is the same as the original product on $\B(\Hi\oplus\Hi)$, and $\LM(X)=\left[\begin{matrix}\B(\Hi)&\Kp(\Hi)\\\B(\Hi)&\Kp(\Hi)^1\end{matrix}\right]$, and accordingly $\Q(X)=\left[\begin{matrix}\B(\Hi)&\B(\Hi)\\\B(\Hi)&\Kp(\Hi)^1\end{matrix}\right]$ which is not a TRO.
\end{example}
\section{The quasi-multiplier space of a dual operator space}\label{section: dual}
In this section, we prove the following theorem. The argument is parallel to that of \cite{B6}~Corollary~3.2~(1).
\begin{theorem}\label{th: dual}If $X$ is an operator space with an operator space predual, then so is $\Q(X)$. Thus by the Banach-Alaoglu Theorem $\Ball(\Q(X))$ is compact in the weak$^*$ topology, and hence by the Krein-Milman Theorem $\Ball(\Q(X))$ is the weak$^*$-closure of the convex hull of the extreme points of $\Ball(\Q(X))$.
\end{theorem}
To prove this, we need several lemmas. Note that if $X$ is an operator space with an operator space predual $X_*$, then $\M_n(X)$ also has an operator space predual which is given by the operator space projective tensor product $\Tr_n\widehat{\otimes}X_*$, where $\Tr_n$ is the set of $n\times n$ trace-class matrices, i.e., $\Tr_n=\M_n(\C)$ as vector spaces, but $\Tr_n$ is given an operator space structure by the identification $\Tr_n\cong\M_n(\C)^*$ with the pairing $<\alpha,\beta>:=\sum_{i,j}\alpha_{i,j}\beta_{i,j},\;\forall\alpha=[\alpha_{i,j}]\in\Tr_n,\;\forall\beta=
[\beta_{i,j}]\in\M_n(\C)$.
\begin{lemma}\label{lm: w*}{\em (\cite{B6}~Lemma~1.6)} If $X$ is a dual operator space, and $x_i$ is a net in $\M_n(X)$, then $x_i\to x\in\M_n(X)$ in the weak$^*$ topology of $\M_n(X)$ if and only if each entry in $x_i$ converges in $X$ in the weak$^*$ topology of $X$ to the corresponding entry in x.
\end{lemma}
\begin{lemma}\label{lm: cpt}{\em (\cite{B6}~Lemma~3.1)} Let $X$ and $Y$ be operator spaces, with $Y$ a dual operator space, and let $T:X\to Y$ be a one-to-one linear mapping. Then the following are equivalent:
\begin{itemize}
\item[(i)]$X$ has an operator space predual such that $T$ is weak$^*$-continuous;
\item[(ii)]$T^{(n)}(\Ball(\M_n(X)))$ is weak$^*$-compact for every positive integer $n$.\footnote{$T^{(n)}:\M_n(X)\to\M_n(Y)$ is defined by $T^{(n)}([x_{i,j}]):=[T(x_{i,j})],\;\forall[x_{i,j}]\in\M_n(X)$.}
\end{itemize}
\end{lemma}
\begin{lemma}\label{lm: matrix qm}Let $X$ be an operator space. Then $\M_n(\Q(X))\cong\Q(\M_n(X))$, completely isometrically.
\end{lemma}
\begin{proof}The assertion easily follows from $\M_n(I(X))\cong I(\M_n(X))$ completely isometric, and the definition of the quasi-multipliers (Definition~\ref{def: qm}).
\end{proof}
{\em Proof of Theorem~\ref{th: dual}.} It suffices to show that the completely contractive one-to-one mapping $\iota:\Q(X)\to\CB(X\stackrel{h}{\otimes}X,X)$ defined by $\iota(z)(x\otimes y):=xzy$ and the completion of their span, satisfies (ii) of Lemma~\ref{lm: cpt}. We need to show that if $\{\varphi^\lambda\}$ is a net in $\Ball(\M_n(\Q(X)))$ converging in the weak$^*$ topology to $\varphi\in\M_n(\CB(X\stackrel{h}{\otimes}X,X))$, then $\varphi\in\Ball(\M_n(\Q(X)))$, where we are identifying $\CB(X\stackrel{h}{\otimes}X,X)$ with $((X\stackrel{h}{\otimes}X)\widehat{\otimes}X_*)^*$ completely isometrically. But by Lemma~\ref{lm: w*} and the canonical identification of Lemma~\ref{lm: matrix qm}, it is enough to show that if $\varphi^\lambda$ is a net in $\Ball(\Q(X))$ converging in the weak$^*$ topology to $\varphi\in\CB(X\stackrel{h}{\otimes}X,X)$, then $\varphi\in\Ball(\Q(X))$. Let $x:=[x_{p,q}],y:=[y_{p,q}],v:=[v_{p,q}],w:=[w_{p,q}]\in\M_m(X),\;1\le p,q\le m$. By Theorem~\ref{th: kaneda}~(iii)$\Rightarrow$(ii), we are done if we have shown that \begin{equation}\label{eq: 1}\left\|\left[\begin{matrix}v_{p,q}&\sum_{k_{p,q}}\varphi(x_{p,q}^{(k_{p,q})},y_{p,q}^{(k_{p,q})})\\0&w_{p,q}
\end{matrix}\right]\right\|\le\left\|\left[\begin{matrix}v_{p,q}\otimes1&\sum_{k_{p,q}}x_{p,q}^{(k_{p,q})}\otimes y_{p,q}^{(k_{p,q})}\\0&1\otimes w_{p,q}\end{matrix}\right]\right\|,\end{equation}where each matrix is $2m\times2m$. However we do know by Theorem~\ref{th: kaneda}~(ii)$\Rightarrow$(iii) that
\begin{equation}\label{eq: 2}\left\|\left[\begin{matrix}v_{p,q}&\sum_{k_{p,q}}\varphi^\lambda(x_{p,q}^{(k_{p,q})},y_{p,q}^{(k_{p,q})})\\0&w_{p,q}
\end{matrix}\right]\right\|\le\left\|\left[\begin{matrix}v_{p,q}\otimes1&\sum_{k_{p,q}}x_{p,q}^{(k_{p,q})}\otimes y_{p,q}^{(k_{p,q})}\\0&1\otimes w_{p,q}\end{matrix}\right]\right\|.
\end{equation}Since $\varphi^\lambda\to\varphi$ in the weak$^*$ topology in $((X\stackrel{h}{\otimes}X)\widehat{\otimes}X_*)^*$, $\varphi^\lambda(x,y)\to\varphi(x,y),\;\forall x,y\in X$ in the weak$^*$ topology in $X$. Indeed, $\varphi^\lambda\stackrel{\operatorname{w}^*}{\to}\varphi$ means that $\forall x,y\in X,\;\forall f\in X_*$, $<\varphi^\lambda,(x\otimes y)\otimes f>\to<\varphi,(x\otimes y)\otimes f>$. But $<\varphi^\lambda,(x\otimes y)\otimes f>=<\varphi^\lambda(x,y),f>$ and $<\varphi,(x\otimes y)\otimes f>=<\varphi(x,y),f>$, thus $<\varphi^\lambda(x,y),f>\to<\varphi(x,y),f>,\;\forall x,y\in X,\;\forall f\in X_*$. Let us denote the matrix of the right hand side of Equation~(\ref{eq: 1}) or (\ref{eq: 2}) by $[\xi_{r,s}]$, and the matrix of the left hand side of Equation~(\ref{eq: 1}) by $[x_{r,s}]$, and the matrix of the left hand side of
Equation~(\ref{eq: 2}) by $[x_{r,s}^\lambda]$. Let $G\in\Ball(\M_{2m}(X)_*)$ which can be identified with $[g_{r,s}]\in \Ball(\Tr_{2m}\widehat{\otimes}X_*)$. Then$$|<[x_{r,s}^\lambda],G>|=\left|\sum_{r,s=1}^{2m}<x_{r,s}^\lambda,g_{r,s}>\right|\le\|[\xi_{r,s}]\|.$$By taking the limit $\lambda\to\infty$, we have that $$|<[x_{r,s}],G>|=\left|\sum_{r,s=1}^{2m}<x_{r,s},g_{r,s},>\right|\le\|[\xi_{r,s}]\|.$$Since $G\in\Ball(\M_{2m}(X)_*)$ is arbitrary, $\|[x_{r,s}]\|\le\|[\xi_{r,s}]\|$, i.e., Inequality~(\ref{eq: 1}) has been shown.\begin{flushright}$\square$
\end{flushright}From the proof above, the following corollary immediately follows.
\begin{corollary}\label{cor:dual}Let $X$ be a dual operator space, and $\{z_i\}\subset\Q(X)$ be a bounded net, and $z\in\Q(X)$. Then $z_i\to z$ in the weak$^*$ topology of $\Q(X)$ if and only if $xz_iy\to xzy,\;\forall x,y\in X$ in the weak$^*$ topology of $X$.
\end{corollary}
\section{Ideal decompositions of a TRO with predual}\label{section:decompositions}
In this section, we prove that any TRO with predual can be decomposed to the direct sum of a two-sided ideal, a left ideal, and a right ideal in some von Neumann algebra. Although this theorem is not directly related to quasi-multipliers, the main tool used to prove it is an extreme point, and the result itself is interesting, so it would be appropriate to present in the present paper. In the special case that a TRO is finite-dimensional, the TRO is decomposed to the direct sum of rectangular matrices, which was essentially first proved by R.~R.~Smith (\cite{S}). We included this result in Appendix with the author's short proof.
\begin{theorem}\label{th: ideal decomposition}Let $X$ be a TRO which is also a dual Banach space. Then $X$ can be
decomposed to the direct sum of TRO's $X_T$, $X_L$, and $X_R$:$$X=X_T\stackrel{\infty}{\oplus}X_L\stackrel{\infty}{\oplus}X_R$$ so that there is a complete isometry $\iota$ from $X$ into a von Neumann algebra in which $\iota(X_T)$, $\iota(X_L)$, and $\iota(X_R)$ are a weak$^*$-closed two-sided, left, and right ideal, respectively, and$$\iota(X)=\iota(X_T)\stackrel{\infty}{\oplus}\iota(X_L)\stackrel{\infty}{\oplus}\iota(X_R).$$
\end{theorem}
\begin{proof}By \cite{EOR}~Theorem~2.6, we may regard $X$ as a weak$^*$-closed subspace of $\B(\K,\Hi)$ for some Hilbert spaces $\Hi$ and $\K$ such that $XX^*X\subset X$. We may assume that $[X\K]=\Hi$ and $[X^*\Hi]=\K$. We also identify $\B(\K,\Hi)$ with the ``(1,2)-corner'' of $\B(\Hi\oplus\K)$, and let $1_{\Hi}\in\B(\Hi\oplus\K)$ and $1_{\K}\in\B(\Hi\oplus\K)$ denote the orthogonal projections on $\Hi$ and $\K$, respectively. Then$$\Li(X):=
\left[
  \begin{array}{cc}
    \overline{XX^*}^{\operatorname{w}^*} & X \\
    X^* & \overline{X^*X}^{\operatorname{w}^*} \\
  \end{array}
\right]
$$is the linking von Neumann algebra, and $1_{\Hi},1_{\K}\in\Li(X)$, and $X=1_{\Hi}\Li(X)1_{\K}$. Since $\Ball(X)$ is weak$^*$-closed in $\B(\K,\Hi)$, there is an extreme point $e\in\Ball(X)$. By Kadison's theorem (Lemma~\ref{lm: kadison}),
\begin{equation}\label{eq:Kadison}(1_{\Hi}-ee^*)X(1_{\K}-e^*e)=\{0\},
\end{equation}and $e$ is a partial isometry. Let $p$ and $q$ be the identities of the von Neumann algebras $\overline{X(1_{\K}-e^*e)X^*}^{\operatorname{w}^*}$ and $\overline{X^*(1_{\Hi}-ee^*)X}^{\operatorname{w}^*}$, respectively. Then by Equation~(\ref{eq:Kadison}), it follows that
\begin{equation}\label{eq:null}pXq=\{0\},
\end{equation}
\begin{equation}\label{eq:commute}p=pee^*=ee^*p=pee^*p,\quad q=qe^*e=e^*eq=qe^*eq,\quad\text{and}
\end{equation}
\begin{equation}\label{eq:central}pxy^*=pxy^*p=xy^*p,\quad qx^*y=qx^*yq=x^*yq,\quad\forall x,y\in X.
\end{equation}Put $q_1:=e^*(1_{\Hi}-p)e(1_{\K}-q)$ and $q_2:=1_{\K}-q-q_1$.

We claim that $q_1$ and $q_2$ are orthogonal projections. Equations~(\ref{eq:central}) noting that $pe\in X$ yield that $q_1^*=(1_{\K}-q)e^*(1_{\Hi}-p)e=e^*e-e^*pe-qe^*e+qe^*pe=e^*e-e^*pe-e^*eq+e^*peq=q_1$ and $q_1^2=e^*(1_{\Hi}-p)e(1_{\K}-q)e^*(1_{\Hi}-p)e(1_{\K}-q)=e^*(1_{\Hi}-p)eq_1^*(1_{\K}-q)=e^*(1_{\Hi}-p)eq_1(1_{\K}-q)=
e^*(1_{\Hi}-p)ee^*(1_{\Hi}-p)e(1_{\K}-q)(1_{\K}-q)=e^*ee^*(1_{\Hi}-p)(1_{\Hi}-p)e(1_{\K}-q)(1_{\K}-q)=e^*(1_{\Hi}-p)e
(1_{\K}-q)=q_1$. Noting that $q_1q=0$, we have that $q_2^2=q_2=q_2^*$. Hence $q_1$ and $q_2$ are orthogonal projections.

To see that
\begin{equation}\label{eq:null2}(1_{\Hi}-p)X(1_{\K}-e^*e)=\{0\},
\end{equation}let $\{u_{\alpha}\}$ be an approximate identity of the $C^*$-algebra $X^*X$. Then for each $x\in X$, $px(1_{\K}-e^*e)u_{\alpha}=x(1_{\K}-e^*e)u_{\alpha}$. Taking the limit $\alpha\to\infty$ yields that $px(1_{\K}-e^*e)=x(1_{\K}-e^*e),\;\forall x\in X$, and hence Equation~(\ref{eq:null2}) holds. Similarly,
\begin{equation}\label{eq:null3}(1_{\Hi}-ee^*)X(1_{\K}-q)=\{0\}
\end{equation}also holds.

Let $x,y\in X$. Then
\begin{center}
\begin{tabular}{rll}
  $q_1x^*y$&$=e^*(1_{\Hi}-p)e(1_{\K}-q)x^*y$&\\
&$=e^*(1_{\Hi}-p)ex^*y(1_{\K}-q)$&\quad by Equation~(\ref{eq:central})\\
&$=e^*ex^*(1_{\Hi}-p)y(1_{\K}-q)$&\quad by Equation~(\ref{eq:central})\\
&$=x^*(1_{\Hi}-p)y(1_{\K}-q)$&\quad by Equation~(\ref{eq:null2})\\
&$=x^*(1_{\Hi}-p)ye^*e(1_{\K}-q)$&\quad by Equation~(\ref{eq:null2})\\
&$=x^*ye^*(1_{\Hi}-p)e(1_{\K}-q)$&\quad by Equation~(\ref{eq:central})\\
&$=x^*yq_1$,&
\end{tabular}
\end{center} and so we have that
\begin{equation}\label{eq:central2}q_1x^*y=x^*yq_1=q_1x^*yq_1,\quad\forall x,y\in X.
\end{equation}

Put $X_T:=Xq_1$, $X_L:=Xq$, and $X_R:=Xq_2$, then these are weak$^*$-closed TRO's, and $X=X_T\oplus X_L\oplus X_R$. Using Equations~(\ref{eq:central})~and~(\ref{eq:central2}) and noting that $q_1$, $q$, and $q_2$ are mutually disjoint, we have that $X_T^*X_L=X_T^*X_R=X_L^*X_T=X_L^*X_R=X_R^*X_T=X_R^*X_L=\{0\}$ and $X^*X=X_T^*X_T\stackrel{\infty}{\oplus}X_L^*X_L\stackrel{\infty}{\oplus}X_R^*X_R$. This proves that $X=X_T\stackrel{\infty}{\oplus}X_L\stackrel{\infty}{\oplus}X_R$.

Define $\iota:X\to\overline{XX^*}^{\operatorname{w}^*}\stackrel{\infty}{\oplus}\overline{X^*X}^{\operatorname{w}^*}$ by$$\iota(x):=(x_T+x_L)e^*\oplus e^*x_R,$$where $x=x_T+x_L+x_R$ is the unique decomposition of $x\in X$ such that $x_T\in X_T$, $x_L\in X_L$, and $x_R\in X_R$. First note that $\iota(X_T)\cap\iota(X_L)=\{0\}$. Indeed, assume that $\iota(x_T)+\iota(x_L)=0$, i.e., $xq_1e^*+xqe^*=0$. Then by multiplying both sides by $e$ on the right and using Equations~(\ref{eq:commute})~and~(\ref{eq:central2}), we obtain that $xe^*eq_1+xq=0$. Multiplying both sides by $q$ on the right noting that $q_1q=0$ yields that $xq=0$, and hence $xq_1e^*=xqe^*=0$, i.e., $\iota(x_T)=\iota(x_L)=0$. Since $\iota(X_T)^*\iota(X_L)=eX_T^*X_Le^*=\{0\}$ and $\iota(X_L)^*\iota(X_T)=eX_L^*X_Te^*=\{0\}$, $(\iota(X_T)\oplus\iota(X_L))^*(\iota(X_T)\oplus\iota(X_L))=\iota(X_T)^*
\iota(X_T)\stackrel{\infty}{\oplus}\iota(X_L)^*\iota(X_L)$ noting that $\iota(X_T)^*\iota(X_T)=q_1X_T^*X_Tq_1$ and $\iota(X_L)^*\iota(X_L)=qX_L^*X_Lq$. Thus $\iota(X)=\iota(X_T)\stackrel{\infty}{\oplus}\iota(X_L)\stackrel{\infty}{\oplus}\iota(X_R)$. To show that $\iota$ is a complete isometry, it suffices to show that each of $\iota_{|{X_T}}$, $\iota_{|{X_L}}$, and $\iota_{|{X_R}}$ is a complete isometry. Since $e^*eq_1=q_1$, $\|\iota(x_T)\|^2=\|\iota(x_T)\iota(x_T)^*\|=\|xq_1e^*eq_1x^*\|=\|xq_1x^*\|=\|xq_1\|^2=\|x_T\|^2$. A similar calculation works at the matrix level, which concludes that $\iota_{|{X_T}}$ is a complete isometry. Similarly, Equation~(\ref{eq:commute}) yields that $\iota_{|{X_L}}$ is a complete isometry. $\|\iota(x_R)\|^2=\|\iota(x_R)^*\iota(x_R)\|=\|q_2x^*ee^*xq_2\|=\|q_2x^*ee^*x(1_{\K}-q-q_1)\|=\|q_2x^*x(1_{\K}-q)\|=
\|q_2x^*x(1_{\K}-q-q_1)\|=\|q_2x^*xq_2\|=\|x_R\|^2$, where we used Equations~(\ref{eq:null3})~and~(\ref{eq:central2}) as well as the fact that $q_2q_1=0$ at the fourth equality, and Equation~(\ref{eq:central2}) together with the fact that $q_2q_1=0$ at the fifth equality. A similar calculation works at the matrix level, which concludes that $\iota_{|{X_R}}$ is a complete isometry.

By \cite{B6}~Lemma~1.5~(3) or \cite{BL}~Theorem~A.2.5~(3) for example, $\iota(X_T)$, $\iota(X_L)$, and $\iota(X_R)$ are weak$^*$-closed. Clearly, $\iota(X_T)$ and $\iota(X_L)$ are left ideals and $\iota(X_R)$ is a right ideal in the von Neumann algebra $\overline{XX^*}^{\operatorname{w}^*}\stackrel{\infty}{\oplus}\overline{X^*X}^{\operatorname{w}^*}$. To see that $\iota(X_T)$ is a right ideal as well, it suffices to show that $\iota(X_T)^*\subset\iota(X_T)$, in which case necessarily $\iota(X_T)^*=\iota(X_T)$. To show this, first note that it follows from Equation~(\ref{eq:null3}) that $q_1x^*=e^*(1_{\Hi}-p)e(1_{\K}-q)x^*=e^*(1_{\Hi}-p)e(1_{\K}-q)x^*ee^*=q_1x^*ee^*,\;\forall x\in X$. Therefore, together with Equations~(\ref{eq:central2}), we have that $\forall x\in X,\;\iota(x_T)^*=eq_1x^*=eq_1x^*ee^*=ex^*eq_1e^*\in Xq_1e^*=\iota(X_T)$.
\end{proof}
\begin{definition}
We call the decomposition $X=X_T\stackrel{\infty}{\oplus}X_L\stackrel{\infty}{\oplus}X_R$ obtained in the proof of Theorem~\ref{th: ideal decomposition} the \emph{\textbf{ideal decomposition}} of the TRO $X$ with predual with respect to the extreme point $e$ of $\Ball(X)$.
\end{definition}
\begin{remark}
\begin{enumerate}
\item The reader should distinguish ideal decompositions from Peirce decompositions in the literature of Jordan triples. In fact, a TRO can be regarded as a Jordan triple with the canonical symmetrization of the triple product. However, an ideal decomposition and a Peirce decomposition give totally different decompositions.
\item It is also possible to define $\iota:X\to\overline{XX^*}^{\operatorname{w}^*}\stackrel{\infty}{\oplus}\overline{X^*X}^{\operatorname{w}^*}$ by $\iota(x):=x_Le^*\oplus e^*(x_R+x_T),\;x\in X$.
\item Simpler expressions for $X_T$ and $X_R$ are $X_T=\{x-px-xq\;;\;x\in X\}$ and $X_R=pX$, respectively, which would be more helpful in understanding what is going on in the decomposition. To see the equivalences of expressions, let $x\in X$. Then using Equations~(\ref{eq:central}),~(\ref{eq:null2}),~and~(\ref{eq:null}), we have that $x_T:=xq_1=xe^*(1_{\Hi}-p)e(1_{\K}-q)=(1_{\Hi}-p)xe^*e(1_{\K}-q)=(1_{\Hi}-p)x(1_{\K}-q)=x-px-xq$. Accordingly, it follows that $x_R:=xq_2=x(1_{\K}-q-q_1)=x(1_{\K}-q)-xq_1=x(1_{\K}-q)-(x-px-xq)=px$.
\item The ideal decomposition highly depends on the extreme point chosen. Indeed, let $X$ be a von Neumann algebra, $u\in X$ be a unitary element, and $w\in X$ be an isometry which is not unitary. Then the ideal decomposition associated with $u$ is just $X=X_T$, while the one associated with $w$ is $X=X_T\stackrel{\infty}{\oplus}X_L$.
\end{enumerate}
\end{remark}
\section*{Appendix: A short proof of Smith's result}The following theorem was first proved by R.~R.~Smith in \cite{S}. However, the author observed it independently as well as Corollary~A.2 in early 2000 when the author was unaware of Smith's result. Since these results are a spacial case of Theorem~\ref{th: ideal decomposition} in the present paper, and also our proof is short enough to understand the essence of the results transparently, it seems worthwhile to present here, although these are already included in \cite{ER2} (Lemma~6.1.7 and Corollary~6.1.8). The key to the extreme shortness of the proof is to note the obvious fact that if a TRO $X$ is finite-dimensional, then so are the $C^*$-algebras $XX^*$ and $X^*X$.
\begin{flushleft}
\textbf{Theorem A.1.} (Smith \cite{S}) {\em If $X$ is a finite-dimensional TRO, then there exist a finite-dimensional $C^*$-algebra $\A$ and an orthogonal projection $p\in\A$ such that $X\cong p\A p^{\perp}$ completely isometrically.}
\end{flushleft}
\begin{proof}Let $X\subset\B(\Hi)$ be a finite-dimensional TRO and $\{x_1,\dots,x_n\}\subset X$ be its base. Then the $C^*$-algebra $XX^*:=\operatorname{span}\{xy^*\;;\;x,y\in X\}$ is equal to the set $\operatorname{span}\{x_ix_j^*\;;\;1\le i,j\le n\}$, and the latter is obviously a finite-dimensional vector space. Similarly, $X^*X:=\operatorname{span}\{x^*y\;;\;x,y\in X\}$ is a finite-dimensional $C^*$-algebra. Let $\Li(X)$ be the linking $C^*$-algebra for $X$, i.e., $\Li(X):=\left[\begin{array}{cc}XX^*&X\\X^*&X^*X\\\end{array}\right]\left(\subset\M_2(\B(\Hi))\right)$ with the canonical $C^*$-algebra structure inherited from $\B(\Hi\oplus\Hi)$. Let $e,f\in\B(\Hi)$ be the identities of the $C^*$-algebras $XX^*$ and $X^*X$, respectively, and let $p:=\left[\begin{array}{cc}e&0\\0&0\\\end{array}\right]\in\Li(X)$. Then $p^{\perp}=\left[\begin{array}{cc}0&0\\0&f\\\end{array}\right]$ and $X\cong p\Li(X)p^{\perp}$ completely isometrically.
\end{proof}
\begin{flushleft}
\textbf{Corollary A.2.} {\em A finite-dimensional TRO is completely isometric to the direct sum of rectangular matrices: $\M_{l_1,k_1}(\C)\stackrel{\infty}{\oplus}\cdots\stackrel{\infty}{\oplus}\M_{l_m,k_m}(\C)$.}
\end{flushleft}
\begin{proof}Let $X$ be a finite-dimensional TRO. By Theorem~A.1, we may assume that $X=p\left(\bigoplus_{i=1}^m\M_{n_i}(\C)\right)p^{\perp}$, where $p$ is an orthogonal projection in $\bigoplus_{i=1}^m\M_{n_i}(\C)$. For each $1\le i\le m$, let us denote by $1_i$ the identity of $\M_{n_i}(\C)$ which is identified with an element of $\bigoplus_{i=1}^m\M_{n_i}(\C)$ in the obvious way, and let $p_i:=p1_i$. Then $X=\bigoplus_{i=1}^mp_i\M_{n_i}(\C)p_i^{\perp}$. By a unitary transform which is a complete isometry, we may assume that $p_i=\operatorname{diag}\{\stackrel{l_i}{\overbrace{1,\dots,1}},0,\dots,0\}$ and $p_i^{\perp}=\operatorname{diag}\{\stackrel{l_i}{\overbrace{0,\dots,0}},1,\dots,1\}$ for each $1\le i\le m$.
\end{proof}

\end{document}